\documentclass[10pt]{amsart}
\usepackage{amsmath,amssymb,amstext,amsthm,amscd}

\theoremstyle{plain}

\newtheorem{theorem}{Theorem}
\newtheorem*{theorem*}{Theorem}

\newtheorem*{claim*}{Claim}
\newtheorem{conjecture}{Conjecture}
\newtheorem*{conjecture*}{Conjecture}
\newtheorem*{question*}{Question}

\theoremstyle{definition}

\newtheorem*{definition*}{Definition}

\theoremstyle{remark}

\newtheorem*{remark*}{Remark}

\newtheorem*{acknowledgement}{Acknowledgement}

\newcommand{\QQ}{\mathbb{Q}}

\newcommand{\Mbar}{\overline{\mathcal{M}}}

\def\({\left(}
\def\){\right)}
\def\<{\langle}
\def\>{\rangle}
\def\d{\partial}

\def\cM{{\mathcal M}}
\def\ocM{{\overline{\cM}}}

\thanks{The second author is partially supported by NSF}

\begin{document}

\title{Tautological equations in genus $2$ via invariance conjectures}

\author{D.~Arcara}

\author{Y.-P.~Lee}

\address{Department of Mathematics, University of Utah,
Salt Lake City, UT 84112-0090, USA}

\email{arcara@math.utah.edu}

\email{yplee@math.utah.edu}

\begin{abstract}
We verify the \emph{Invariance Conjectures} of tautological equations 
\cite{ypL2} in genus two.
In particular, a uniform derivation of all known genus two equations is given.
\end{abstract}

\maketitle

\setcounter{section}{-1}


\section{Introduction}
The purpose of this paper is to verify the genus two case of 
\emph{Invariance Conjectures} of tautological equations proposed in 
\cite{ypL2}.
In particular, applying Theorem~5 in \cite{ypL3} 
(i.e.~Conjecture~1 in \cite{ypL2}) and E.~Getzler's Hodge numbers 
calculations \cite{eG2}, we are able to give
a uniform derivation of all known genus two tautological equations:
Mumford--Getzler's equation, Getzler's equation \cite{eG2} and
Belorousski--Pandharipande's equation \cite{BP}.
This, combined with \cite{GL} in genus one and \cite{AL} in genus three, 
shows that this method generates and proves all known tautological equations.

\subsection{Review of the Conjectures}
Some ingredients in \cite{ypL2} are needed in our calculations.

\begin{itemize}
\item The elements in the $\QQ$-tautological algebra $R^k(\Mbar_{g,n})$ 
are the tautological classes of codimension $k$ in $\Mbar_{g,n}$.
The notation $(g,n,k)$ will be used to denote the triple 
(genus, number of marked points, codimension).
\item Any tautological class can be represented by a linear combination
of decorated (and labeled) graph $\Gamma$. See Section~1 in \cite{ypL2}.
We use a special notation of graphs, called \emph{gwi}s, which 
are explained in Section~3.1 in \cite{ypL2}.
\item The operators
\[ 
 \mathfrak{r}_l: R^k(\Mbar_{g,n}) \to R^{k+l-1}(\Mbar^{\bullet}_{g-1,n+2}),
 \quad l=1,2,\ldots
\]
are defined in Definition~2 of \cite{ypL2} as operations on decorated graphs.
Note that the image lies in the moduli of possibly \emph{disconnected} curves,
with at most two connected components.
\end{itemize}

Let $\sum_i c_i \Gamma_i =0$
be a tautological equation in codimension $k$ strata in $\overline{M}_{g,n}$.
The Invariance Conjectures in \cite{ypL2} are:

\setcounter{conjecture}{+1}

\begin{theorem} {\rm (\cite{ypL3}, originally Conjecture~1)}
For all $l$
\begin{equation} \label{e:1}
 \mathfrak{r}_l (\sum_i c_i \Gamma_i) =0.
\end{equation}
\end{theorem}

Let $E=\sum_i c_i \Gamma_i$ be a given linear combination of codimension $k$
tautological strata in $\Mbar_{g,n}$ and $k < 3g-3+n$, with 
$c_i$ \emph{unknown variables}.

\begin{conjecture}
If $\mathfrak{r}_l (E) =0$ for all $l$, then $E=0$ is a tautological equation.
\end{conjecture}

\begin{conjecture}
Conjecture~2 will produce all tautological equations inductively.
\end{conjecture}

\subsection{The algorithm of finding tautological equations}

A general algorithm of finding the tautological equations, based on
Conjectures~2 and 3, is explained in \cite{ypL2} Section~2.
Since these remain conjectural, one possible alternative to the general scheme
is to
\begin{itemize}
\item Calculate the rank of $R^k(\Mbar_{g,n})$ to 
see if there is any new equation.
\item If there is one, then apply invariance condition equation \eqref{e:1} 
to obtain the coefficients of the equation.
\end{itemize}

Since Theorem~1 gives a \emph{necessary} condition, this proceedure gives
a proof of the generated tautological equation.

\subsection{Main results}

\begin{theorem}
Invariance Conjectures hold for $(g,n,k)=(2,1,2),(2,2,2),(2,3,2)$.
In particular, a uniform derivation of all known genus two tautological 
equations is given by invariance condition \eqref{e:1}.
\end{theorem}

\begin{remark*}
As explained in \cite{ypL2} and \cite{ypL3},
our calculation in terms of gwis can be translated literally into one for
any (axiomatic) Gromov--Witten theories.
Therefore, it completes (the write-up of) a proof of the genus two case of
Virasoro conjecture in the semisimple case and of
Witten's conjecture (on spin curves and Gelfand--Dickey hierarchies).
\end{remark*}

\begin{acknowledgement}
We wish to thank A.~Bertram, E.~Getzler, A.~Givental, R.~Pandharipande, and
R.~Vakil for many useful discussions. 
The final stage of this work was done during
the second author's visit to NCTS, whose hospitality is greatly appreciated.
\end{acknowledgement}

\section{Mumford--Getzler's equation in $\ocM_{2,1}$}

In all calculations below, we will emply the ``gwi'' notations for decorated
graphs. It is explained in \cite{ypL2} that gwis are equivalent to the
decorated graphs, or a tautological class. The notations are obviously
inspired by Gromov--Witten invariants.

\subsection{Tautological classes of $R^2(\ocM_{2,1})$}

There are $8$ boundary strata of codimension $\leq 2$ in $\ocM_{2,1}$:
$1$ stratum in codimension $0$, $2$ strata in codimension $1$, and $5$ strata
in codimension $2$.
If we insert $\psi$ classes, the $2$ boundary strata in codimension $1$ 
produce $5$ different tautological elements in codimension $2$.
Note that the $\kappa$-classes can be expressed in terms of boundary and
$\psi$-classes in genus two. So the only decoration one would need is the
$\psi$-classes.

Here is a list of all the $11$ strata with $\psi$ classes in codimension $2$:
\[
 \begin{split}
 &\< \d^x_2 \>_2, \ \< \d^x \d^\mu_1 \d^\mu \>_1, 
 \< \d^x_1 \d^\mu \d^\mu \>_1, \
 \< \d^x \d^\mu_1 \>_1 \< \d^\mu \>_1, \
 \< \d^x \d^\mu \>_1 \< \d^\mu_1 \>_1, \
 \< \d^x_1 \d^\mu \>_1 \< \d^\mu \>_1, \\
 &\< \d^\mu \>_1 \< \d^x \d^\mu \d^\nu \d^\nu \>, \
 \< \d^x \d^\mu \>_1 \< \d^\mu \d^\nu \d^\nu \>, \
 \< \d^\mu \d^\nu \>_1 \< \d^x \d^\mu \d^\nu \>, \
 \< \d^\mu \>_1 \< \d^\nu \>_1 \< \d^x \d^\mu \d^\nu \>, \\
 &\< \d^x \d^\mu \d^\mu \d^\nu \d^\nu \>.
\end{split}
\]

The $5$ strata with $\psi$-classes can be written in terms of the $5$ strata
without $\psi$-classes using TRR's, and therefore we only have $6$ terms which
could be independent.
A general element can be written as
\[
 \begin{split}
 E=&c_1  \< \d^x_2 \>_2 + c_2 \< \d^\mu \>_1 \< \d^x \d^\mu \d^\nu \d^\nu \>
 + c_3 \< \d^x \d^\mu \>_1 \< \d^\mu \d^\nu \d^\nu \> \\
 + &c_4 \< \d^\mu \d^\nu \>_1 \< \d^x \d^\mu \d^\nu \>
 + c_5 \< \d^\mu \>_1 \< \d^\nu \>_1 \< \d^x \d^\mu \d^\nu \>
 + c_6 \< \d^x \d^\mu \d^\mu \d^\nu \d^\nu \>.
 \end{split}
\]

\subsection{Calculating $\mathfrak{r}_1(E)$}

\emph{Throughout this paper, the labelings $i,j$ are assumed to be 
symmetrized for $l$ odd, and anti-symmetrized for $l$ even.}

\begin{eqnarray*}
\< \d^x_2 \>_2 & \mapsto &
- \frac{1}{2} \< \d^\mu \>_1 \< \d^x \d^\mu \d^\nu \> \< \d^i \d^j \d^\nu \>
- \frac{1}{48} \< \d^x \d^\mu \d^\mu \d^\nu \> \< \d^i \d^j \d^\nu \>
\\ & & 
- \frac{1}{24}  \< \d^x \d^\mu \d^\nu \> \< \d^i \d^j \d^\mu \d^\nu \> 
- \frac{1}{24} \< \d^j \>_1 \< \d^x \d^\mu \d^\nu \> \< \d^i \d^\mu \d^\nu \>
\end{eqnarray*}
\begin{eqnarray*}
\< \d^\mu \>_1 \< \d^x \d^\mu \d^\nu \d^\nu \> & \mapsto &
\frac{1}{24} \< \d^i \d^\mu \d^\mu \> \< \d^x \d^j \d^\nu \d^\nu \>
+ \< \d^j \>_1 \< \d^i \d^\mu \d^\nu \> \< \d^x \d^\mu \d^\nu \>
\\ & & 
- \frac{1}{2} \< \d^i \d^j \d^\mu \> \< \d^x \d^\mu \d^\nu \d^\nu \> 
- \< \d^\mu \>_1 \< \d^x \d^i \d^\mu \> \< \d^j \d^\nu \d^\nu \>
\end{eqnarray*}
\begin{eqnarray*}
\< \d^x \d^\mu \>_1 \< \d^\mu \d^\nu \d^\nu \> & \mapsto &
\< \d^\mu \>_1 \< \d^i \d^x \d^\mu \> \< \d^j \d^\nu \d^\nu \>
+ \frac{1}{24} \< \d^i \d^x \d^\mu \d^\mu \> \< \d^j \d^\nu \d^\nu \>
\\ & & 
- \frac{1}{2} \< \d^i \d^j \d^x \d^\mu \> \< \d^\mu \d^\nu \d^\nu \> 
- \< \d^i \>_1 \< \d^j \d^x \d^\mu \> \< \d^\mu \d^\nu \d^\nu \>
\end{eqnarray*}
\begin{eqnarray*}
\< \d^\mu \d^\nu \>_1 \< \d^x \d^\mu \d^\nu \> & \mapsto &
2 \< \d^\mu \>_1 \< \d^i \d^\mu \d^\nu \> \< \d^x \d^j \d^\nu \>
+ \frac{1}{12} \< \d^i \d^\mu \d^\mu \d^\nu \> \< \d^x \d^j \d^\nu \>
\\ & & 
- \frac{1}{2} \< \d^i \d^j \d^\mu \d^\nu \> \< \d^x \d^\mu \d^\nu \> 
- \< \d^i \>_1 \< \d^j \d^\mu \d^\nu \> \< \d^x \d^\mu \d^\nu \> 
\end{eqnarray*}
\begin{eqnarray*}
\< \d^\mu \>_1 \< \d^\nu \>_1 \< \d^x \d^\mu \d^\nu \> & \mapsto &
\frac{1}{12} \< \d^\nu \>_1 \< \d^i \d^\mu \d^\mu \> \< \d^x \d^j \d^\nu \>
- \< \d^\nu \>_1 \< \d^i \d^j \d^\mu \> \< \d^x \d^\mu \d^\nu \>
\end{eqnarray*}
\begin{eqnarray*}
\< \d^x \d^\mu \d^\mu \d^\nu \d^\nu \> & \mapsto &
4 \< \d^i \d^\mu \d^\nu \> \< \d^x \d^j \d^\mu \d^\nu \>
- 2 \< \d^i \d^x \d^\mu \d^\mu \> \< \d^j \d^\nu \d^\nu \> \\
\end{eqnarray*}

%

\subsection{Setting $\mathfrak{r}_1(E)=0$}

Now we will pick a basis, and set its coordinates to zero.

\begin{eqnarray}\label{eq1-1}
 \< \d^\mu \>_1 \< \d^x \d^\mu \d^\nu \> \< \d^i \d^j \d^\nu \>: \quad
- \frac{1}{2} c_1 + 2 c_4 - c_5 = 0.
\end{eqnarray}
\begin{eqnarray}\label{eq1-2}
 \< \d^\mu \>_1 \< \d^x \d^i \d^\mu \> \< \d^j \d^\nu \d^\nu \>: \quad
 - c_2 + c_3 + \frac{1}{12} c_5 = 0.
\end{eqnarray}
\begin{eqnarray}\label{eq1-3}
 \< \d^i \>_1 \< \d^j \d^x \d^\mu \> \< \d^\mu \d^\nu \d^\nu \>: \quad
 - \frac{1}{24} c_1 + c_2 - c_3 - c_4 = 0.
\end{eqnarray}
\begin{eqnarray}\label{eq1-4}
 \< \d^x \d^i \d^\mu \d^\mu \> \< \d^j \d^\nu \d^\nu \>: \quad
 \frac{1}{24} c_2 + \frac{1}{24} c_3 - 2 c_6 = 0.
\end{eqnarray}

The remaining terms are related to each via WDVV as follows:
\[
 \begin{split}
 &\< \d^x \d^i \d^j \d^\mu \> \< \d^\mu \d^\nu \d^\nu \> \\
 = & \< \d^x \d^\mu \d^\mu \d^\nu \> \< \d^i \d^j \d^\nu \>
 +2 \< \d^x \d^\mu \d^\nu \> \< \d^i \d^j \d^\mu \d^\nu \> 
 -2 \< \d^x \d^i \d^\mu \> \< \d^j \d^\mu \d^\nu \d^\nu \>, \\
 &\< \d^x \d^i \d^\mu \d^\nu \> \< \d^j \d^\mu \d^\nu \> \\
 = & \< \d^x \d^\mu \d^\mu \d^\nu \> \< \d^i \d^j \d^\nu \>
 + \< \d^x \d^\mu \d^\nu \> \< \d^i \d^j \d^\mu \d^\nu \> 
 - \< \d^x \d^i \d^\mu \> \< \d^j \d^\mu \d^\nu \d^\nu \>.
 \end{split}
\]
Therefore, among the $5$ vectors, only $3$ of them are linearly independent.

\begin{eqnarray}\label{eq1-5}
 \< \d^x \d^\mu \d^\mu \d^\nu \> \< \d^i \d^j \d^\nu \> : \quad
 - \frac{1}{48} c_1 - \frac{1}{2} c_2 - \frac{1}{2} c_3 + 4 c_6 = 0.
\end{eqnarray}
\begin{eqnarray}\label{eq1-6}
 \< \d^x \d^\mu \d^\nu \> \< \d^i \d^j \d^\mu \d^\nu \>: \quad
 - \frac{1}{24} c_1 - \frac{1}{2} c_4 - c_3 + 4 c_6 = 0.
\end{eqnarray}
\begin{eqnarray}\label{eq1-7}
 \< \d^x \d^i \d^\mu \> \< \d^j \d^\mu \d^\nu \d^\nu\>: \quad
 c_3 + \frac{1}{12} c_4 - 4 c_6 = 0.
\end{eqnarray}

The system of equations (\ref{eq1-1}), (\ref{eq1-2}), (\ref{eq1-3}), (\ref{eq1-4}),
(\ref{eq1-5}), (\ref{eq1-6}), and (\ref{eq1-7}) has a unique solution
(up to scaling)
\[
c_2 = - \frac{13}{240} c_1,\
c_3 = \frac{1}{240} c_1,\
c_4 = - \frac{1}{10} c_1,\
c_5 = - \frac{7}{10} c_1,\
c_6 = - \frac{1}{960} c_1
\]

We therefore obtain that, if we let $c_1=-1$,
\begin{eqnarray*}
& &  - \< \d^x_2 \>_2 
+ \frac{13}{240} \< \d^\mu \>_1 \< \d^x \d^\mu \d^\nu \d^\nu \>
- \frac{1}{240} \< \d^x \d^\mu \>_1 \< \d^\mu \d^\nu \d^\nu \> \\
& & + \frac{1}{10} \< \d^\mu \d^\nu \>_1 \< \d^x \d^\mu \d^\nu \>
+ \frac{7}{10} \< \d^\mu \>_1 \< \d^\nu \>_1 \< \d^x \d^\mu \d^\nu \>
+ \frac{1}{960} \< \d^x \d^\mu \d^\mu \d^\nu \d^\nu \> = 0,
\end{eqnarray*}
which is Mumford--Getzler's equation.

\subsection{Checking $\mathfrak{r}_2(E)=0$}

Let us now calculate $\mathfrak{r}_2(E)$.

\[
\begin{split}
\< \d^x_2 \>_2 &\mapsto  
- \frac{1}{576} \< \d^i \d^\mu \d^\mu \> \< \d^x \d^\alpha \d^j \> \< \d^\alpha \d^\nu \d^\nu \> \\
\< \d^\mu \d^\mu \d^\nu \> \< \d^x \d^\nu \>_1 
&\mapsto 
 \frac{1}{24} \< \d^j \d^\mu \d^\mu \> \< \d^i \d^\alpha \d^\nu \> \< \d^x \d^\alpha \d^\nu \> 
- \frac{1}{24} \< \d^j \d^x \d^\mu \> \< \d^i \d^\nu \d^\nu \> \< \d^\alpha \d^\alpha \d^\mu \> \\
\< \d^\mu \d^\nu \>_1 \< \d^x \d^\mu \d^\nu \> 
&\mapsto 
- \frac{1}{24} \< \d^j \d^\nu \d^\mu \> \< \d^i \d^\alpha \d^\alpha  \> \< \d^x \d^\nu \d^\mu \> \\
\< \d^x \d^\mu \d^\mu \d^\nu \d^\nu \> 
& \mapsto 
- \< \d^j \d^\mu \d^\mu \> \< \d^i \d^\alpha \d^\nu \> \< \d^x \d^\alpha \d^\nu \> 
- \< \d^j \d^\mu \d^\mu \> \< \d^i \d^\alpha \d^\nu \> \< \d^x \d^\alpha \d^\nu \> .
\end{split}
\]
The other graphs all have $\mathfrak{r}_2(\Gamma)=0.$

Therefore, $\mathfrak{r}_2(E)=0$ as
$$ \frac{1}{576} c_1 + \frac{1}{24} c_{3} + \frac{1}{24} c_{3} + \frac{1}{24} c_{4} - c_{6} - c_{6} = 0. $$

\section{Getzler's equation in $\ocM_{2,2}$}

\subsection{Tautological classes in $\ocM_{2,2}$ of codimension $2$}
\label{strata}

A general linear combination of codimension $2$ tautological classes in 
$\ocM_{2,2}$ is, after removing the linearly dependent classes from 
the induced equations (TRR's and Mumford--Getzler's),

\[
\begin{split}
E=
& c_1 \< \d^x_1 \d^y_1 \>_2
+ c_2 \< \d^\mu_1 \>_2 \< \d^x \d^y \d^\mu \>
+ c_3 \< \d^x \d^y \d^\mu \d^\mu \d^\nu \d^\nu \>
+ c_4 \< \d^\mu \d^\mu \d^\nu \>_1 \< \d^x \d^y \d^\nu \> \\
+ &c_5 \< \d^x \d^y \d^\mu \d^\mu \d^\nu \> \< \d^\nu \>_1
+ c_6 \< \d^x \d^\mu \d^\mu \d^\nu \> \< \d^y \d^\nu \>_1
+ c_7 \< \d^y \d^\mu \d^\mu \d^\nu \> \< \d^x \d^\nu \>_1 \\
+ &c_8 \< \d^\mu \d^\mu \d^\nu \> \< \d^x \d^y \d^\nu \>_1
+ c_9 \< \d^x \d^\mu \d^\nu \>_1 \< \d^y \d^\mu \d^\nu \>
+ c_{10} \< \d^y \d^\mu \d^\nu \>_1 \< \d^x \d^\mu \d^\nu \> \\
+ &c_{11} \< \d^\mu \d^\nu \>_1 \< \d^x \d^y \d^\mu \d^\nu \>
+ c_{12} \< \d^\mu \d^\nu \>_1 \< \d^x \d^y \d^\mu \> \< \d^\nu \>_1 
+ c_{13} \< \d^x \d^y \d^\mu \d^\nu \> \< \d^\mu \>_1 \< \d^\nu \>_1 \\
+ &c_{14} \< \d^x \d^\mu \d^\nu \> \< \d^\mu \>_1 \< \d^y \d^\nu \>_1 
+ c_{15} \< \d^y \d^\mu \d^\nu \> \< \d^\mu \>_1 \< \d^x \d^\nu \>_1.
\end{split}
\]

\subsection{Setting $\mathfrak{r}_1 (E) =0$} 
The routine calculation of $\mathfrak{r}_1 (E)$ is omitted.
Again, a basis will be chosen and the components set to zero.

\begin{eqnarray}\label{eq2-1}
 \< \d^x \d^i \>_1 \< \d^y \d^\mu \d^\nu \> \< \d^\nu \d^j \d^\mu \>: \quad
c_7 - c_8 - c_9 = 0.
\end{eqnarray}
\begin{eqnarray}\label{eq2-2}
 \< \d^y \d^i \>_1 \< \d^x \d^\mu \d^\nu \> \< \d^\nu \d^j \d^\mu \>: \quad
 c_6 - c_8 - c_{10} = 0.
\end{eqnarray}
\begin{eqnarray}\label{eq2-3}
 \< \d^x \d^\mu \>_1 \< \d^y \d^i \d^\mu \> \< \d^j \d^\nu \d^\nu \>: \quad
 - c_7 + c_8 + \frac{1}{24} c_{15} = 0.
\end{eqnarray}
\begin{eqnarray}\label{eq2-4}
\< \d^y \d^\mu \>_1 \< \d^x \d^i \d^\mu \> \< \d^j \d^\nu \d^\nu \>: \quad
-c_6 + c_8 + \frac{1}{24} c_{14} = 0.
\end{eqnarray}
\begin{eqnarray}\label{eq2-5}
\< \d^i \d^\mu \>_1 \< \d^x \d^y \d^j \> \< \d^\mu \d^\nu \d^\nu \>: \quad
 - \frac{1}{240} c_2 - c_8 = 0.
\end{eqnarray}
\begin{eqnarray}\label{eq2-6}
\< \d^i \d^\mu \>_1 \< \d^x \d^y \d^\mu \> \< \d^j \d^\nu \d^\nu \>: \quad
- c_4 + \frac{1}{24} c_{12} = 0.
\end{eqnarray}
\begin{eqnarray}\label{eq2-7}
\< \d^i \>_1 \< \d^\mu \>_1 \< \d^x \d^y \d^\nu \> \< \d^\nu \d^j \d^\mu \>: \quad
- 2 c_1 - c_2 + 2 c_{13} - c_{14} - c_{15} = 0.
\end{eqnarray}
\begin{eqnarray}\label{eq2-8}
\< \d^\mu \d^\nu \>_1 \< \d^x \d^y \d^i \> \< \d^j \d^\mu \d^\nu \>: \quad
\frac{1}{10} c_2 + c_4 + c_4 - c_{11} = 0.
\end{eqnarray}
\begin{eqnarray}\label{eq2-9}
\< \d^\mu \>_1 \< \d^\nu \>_1 \< \d^x \d^i \d^\mu \> \< \d^y \d^j \d^\nu \>: \quad
- c_1 - 2 c_{13} + c_{14} + c_{15} = 0.
\end{eqnarray}
\begin{eqnarray}\label{eq2-10}
\< \d^\mu \>_1 \< \d^\nu \>_1 \< \d^x \d^y \d^i \> \< \d^j \d^\mu \d^\nu \>: \quad
\frac{7}{10} c_2 + c_{12} - c_{13} = 0.
\end{eqnarray}
\begin{eqnarray}\label{eq2-11}
\< \d^i \d^\mu \d^\mu \d^\nu \> \< \d^x \d^y \d^j \> \< \d^\nu \>_1: \quad
\frac{13}{240} c_2 + c_4 - c_5 + \frac{1}{24} c_{12} = 0.
\end{eqnarray}
\begin{eqnarray}\label{eq2-12}
\< \d^i \d^x \d^\mu \d^\mu \> \< \d^j \d^y \d^\nu \> \< \d^\nu \>_1: \quad
- \frac{1}{24} c_1 - c_5 + c_6 + \frac{1}{24} c_{15} = 0.
\end{eqnarray}
\begin{eqnarray}\label{eq2-13}
\< \d^i \d^y \d^\mu \d^\mu \> \< \d^j \d^x \d^\nu \> \< \d^\nu \>_1: \quad
- \frac{1}{24} c_1 - c_5 + c_7 + \frac{1}{24} c_{14} = 0.
\end{eqnarray}
\begin{eqnarray}\label{eq2-14}
\< \d^i \d^\mu \d^\mu \> \< \d^x \d^y \d^j \d^\nu \> \< \d^\nu \>_1: \quad
- c_5 + c_8 + \frac{1}{12} c_{13} = 0.
\end{eqnarray}
\begin{eqnarray}\label{eq2-19}
\< \d^x \d^y \d^i \> \< \d^j \d^\mu \d^\mu \d^\nu \d^\nu \>: \quad
\frac{1}{960} c_2 - c_3 + \frac{1}{24} c_4 = 0.
\end{eqnarray}
\begin{eqnarray}\label{eq2-20}
 \< \d^x \d^i \d^\mu \d^\mu \> \< \d^y \d^j \d^\nu \d^\nu \> : \quad
- \frac{1}{576} c_1 - 2 c_3 + \frac{1}{24} c_6 + \frac{1}{24} c_7 = 0.
\end{eqnarray}
\begin{eqnarray}\label{eq2-21}
\< \d^x \d^y \d^i \d^\mu \d^\mu \> \< \d^j \d^\nu \d^\nu \>: \quad
- 2 c_3 + \frac{1}{24} c_5 + \frac{1}{24} c_8 = 0.
\end{eqnarray}

The $7$ vectors 
\[
 \begin{split}
  &\< \d^x \d^i \d^\mu \> \< \d^y \d^\mu \d^\nu \d^\nu \> \< \d^j \>_1,\ 
 \< \d^x \d^i \d^\mu \d^\nu \> \< \d^y \d^\mu \d^\nu \> \< \d^j \>_1,\ 
 \< \d^x \d^y \d^\mu \> \< \d^i \d^\mu \d^\nu \d^\nu \> \< \d^j \>_1,\\ 
 &\< \d^x \d^y \d^\mu \d^\nu \> \< \d^i \d^\mu \d^\nu \> \< \d^j \>_1,\ 
 \< \d^x \d^y \d^i \d^\mu \> \< \d^\mu \d^\nu \d^\nu \> \< \d^j \>_1,\ 
 \< \d^x \d^\mu \d^\nu \> \< \d^y \d^i \d^\mu \d^\nu \> \< \d^j \>_1,\\ 
 &\< \d^x \d^\mu \d^\mu \d^\nu \> \< \d^y \d^i \d^\nu \> \< \d^j \>_1
\end{split}
\]
are related by WDVV equations.
There are $4$ independent vectors.
\begin{equation}\label{eq2-15}
\< \d^x \d^i \d^\mu \> \< \d^y \d^\mu \d^\nu \d^\nu \> \< \d^j \>_1: \quad
- \frac{1}{12} c_1 - c_7 + 2 c_5 - c_{11} - c_6 = 0.
\end{equation}
\begin{equation}\label{eq2-16}
\< \d^x \d^i \d^\mu \d^\nu \> \< \d^y \d^\mu \d^\nu \> \< \d^j \>_1: \quad
- \frac{1}{12} c_1 - c_9 + 2 c_5 - c_{11} + c_{10} - 2 c_6 = 0.
\end{equation}
\begin{equation}\label{eq2-17}
\< \d^x \d^y \d^\mu \> \< \d^i \d^\mu \d^\nu \d^\nu \> \< \d^j \>_1: \quad
-\frac{1}{24} c_2 - c_4 - c_5 + \frac{1}{24}c_{12} + c_{11} - c_{10} + c_6 =0.
\end{equation}
\begin{equation}\label{eq2-18}
\< \d^x \d^y \d^i \d^\mu \> \< \d^\mu \d^\nu \d^\nu \> \< \d^j \>_1: \quad
 - c_8 - c_{10} + c_6 = 0.
\end{equation}

All of the other remaining terms are related to each other
via WDVV, TRR's and Getzler's genus one equation.
After applying the above equations, one can write them in terms of a basis.

\begin{equation}\label{eq12_xi}
\< \d^x \d^i \d^\mu \> \< \d^y \d^j \d^\nu \> \< \d^\mu \d^\nu \>_1: \quad 
c_1 + c_2 + 20 c_4 - 24 c_5 + 24 c_7 + 2 c_9 + 26 c_{10} - 2 c_{11} = 0.
\end{equation}
\begin{equation}\label{eq12_xy}
\< \d^x \d^y \d^\mu \> \< \d^i \d^j \d^\nu \> \< \d^\mu \d^\nu \>_1: \quad
c_1 + \frac{1}{2} c_2 + 12 c_4 - 12 c_5 + 12 c_7 + 12 c_{10} - \frac{1}{2} c_{12} = 0.
\end{equation}
\begin{equation}\label{eqx1}
\< \d^y \d^i \d^\mu \> \< \d^j \d^\mu \d^\nu \> \< \d^x \d^\nu \>_1: \quad
- c_1 - \frac{1}{2} c_2 - 10 c_4 + 12 c_5 - 12 c_7 + 2 c_9 - 12 c_{10} - \frac{1}{2} c_{15} = 0.
\end{equation}
\begin{equation}\label{eqy1}
\< \d^x \d^i \d^\mu \> \< \d^j \d^\mu \d^\nu \> \< \d^y \d^\nu \>_1: \quad
- c_1 - \frac{1}{2} c_2 - 10 c_4 + 12 c_5 - 12 c_7 - 10 c_{10} - \frac{1}{2} c_{14} = 0.
\end{equation}
\begin{equation}\label{eqi1}
\< \d^x \d^y \d^\mu \> \< \d^i \d^\mu \d^\nu \> \< \d^j \d^\nu \>_1: \quad
- 4 c_1 - 2 c_2 - 20 c_4 + 24 c_5 - 24 c_7 - 2 c_9 - 26 c_{10} + 2 c_{11} = 0.
\end{equation}
\begin{equation}\label{eq1_xj2}
\< \d^x \d^i \d^\mu \> \< \d^y \d^j \d^\mu \d^\nu \> \< \d^\nu \>_1: \quad
- 3 c_1 - c_2 - 8 c_4 + 12 c_5 - 12 c_7 + 2 c_9 - 10 c_{10} - \frac{1}{2} c_{12} - c_{13} = 0.
\end{equation}
\begin{equation}\label{eq1_x12}
\< \d^y \d^i \d^j \d^\mu \> \< \d^x \d^\mu \d^\nu \> \< \d^\nu \>_1: \quad
\frac{3}{2} c_1 + \frac{1}{2} c_2 + 10 c_4 - 12 c_5 + 12 c_7 - 2 c_9 + 12 c_{10} + c_{13} - \frac{1}{2} c_{14} = 0.
\end{equation}
\begin{equation}\label{eq1_y12}
 \< \d^y \d^i \d^\mu \> \< \d^x \d^j \d^\mu \d^\nu \> \< \d^\nu \>_1: \quad
- \frac{3}{2} c_1 - \frac{1}{2} c_2 + 2 c_4 + 2 c_9 - \frac{1}{2} c_{12} - \frac{1}{2} c_{15} = 0.
\end{equation}
\begin{equation}\label{eq1_i12}
\< \d^x \d^y \d^i \d^\mu \> \< \d^j \d^\mu \d^\nu \> \< \d^\nu \>_1: \quad
c_1 + c_2 + 8 c_4 - 12 c_5 + 12 c_7 + 12 c_{10} + 2 c_{11} + \frac{1}{2} c_{12} - c_{13} = 0.
\end{equation}
\begin{equation}\label{eqxy21}
\< \d^x \d^y \d^\nu \d^\mu \> \< \d^\mu \d^i \d^j \d^\nu \>: \quad
\frac{1}{24} c_2 + \frac{5}{6} c_4 - \frac{5}{12} c_{11} = 0.
\end{equation}
\begin{equation}\label{eqxi21}
\< \d^x \d^i \d^\nu \d^\mu \> \< \d^\mu \d^y \d^j \d^\nu \>: \quad
- \frac{1}{8} c_1 - \frac{1}{24} c_2 + 4 c_3 - \frac{5}{6} c_4 - c_7 - c_{10} + \frac{1}{12} c_{11} = 0.
\end{equation}
\begin{equation}\label{eqx221}
\< \d^x \d^\nu \d^\nu \d^\mu \> \< \d^\mu \d^y \d^i \d^j \>: \quad
- \frac{1}{24} c_1 - \frac{1}{24} c_2 - \frac{5}{6} c_4 + c_5 - \frac{1}{2} c_6 - \frac{1}{2} c_7 - \frac{1}{12} c_{11} = 0.
\end{equation}
\begin{equation}\label{eqxi1}
\< \d^x \d^i \d^\mu \> \< \d^\mu \d^y \d^j \d^\nu \d^\nu \>: \quad
- \frac{1}{16} c_1 - \frac{1}{48} c_2 - \frac{5}{12} c_4 - \frac{5}{12} c_{10} = 0.
\end{equation}
\begin{equation}\label{eqyi1}
\< \d^y \d^i \d^\mu \> \< \d^\mu \d^x \d^j \d^\nu \d^\nu \>: \quad
- \frac{5}{48} c_1 - \frac{1}{16} c_2 - \frac{5}{4} c_4 + c_5 - c_7 + \frac{1}{12} c_9 - \frac{1}{2} c_{10} - \frac{1}{12} c_{11} = 0.
\end{equation}
\begin{equation}\label{eqy21}
\< \d^y \d^\nu \d^\mu \> \< \d^\mu \d^x \d^i \d^j \d^\nu \>: \quad
\frac{1}{24} c_1 + \frac{1}{24} c_2 + \frac{5}{6} c_4 - c_5 + c_7 - \frac{1}{2} c_9 + \frac{1}{2} c_{10} + \frac{1}{12} c_{11} = 0.
\end{equation}
\begin{equation}\label{eqi21}
\< \d^i \d^\nu \d^\mu \> \< \d^\mu \d^x \d^y \d^j \d^\nu \>: \quad
\frac{1}{12} c_1 + \frac{1}{12} c_2 + 4 c_3 + \frac{5}{3} c_4 - 2 c_5 + c_7 + c_{10} + \frac{1}{12} c_{11} = 0.
\end{equation}
\begin{equation}\label{eq221}
\< \d^\nu \d^\nu \d^\mu \> \< \d^\mu \d^x \d^y \d^i \d^j \>: \quad
- \frac{1}{24} c_1 - \frac{1}{24} c_2 - \frac{5}{6} c_4 + c_5 - \frac{1}{2} c_7 - \frac{1}{2} c_8 - \frac{1}{2} c_{10} - \frac{1}{12} c_{11} = 0.
\end{equation}

Solving equations (\ref{eq2-1})-(\ref{eq221}), we can write all of the
coefficients in terms of $c_1$:
\[
\begin{split}
&c_2 = - 3 c_1, \
c_3 = - \frac{1}{576} c_1,  \
c_4 = \frac{1}{30} c_1,  \
c_5 = - \frac{23}{240} c_1,  \
c_6 = - \frac{1}{48} c_1,  \
c_7 = - \frac{1}{48} c_1,  \\
&c_8 = \frac{1}{80} c_1,  \
c_9 = - \frac{1}{30} c_1,  \
c_{10} = - \frac{1}{30} c_1,  \
c_{11} = - \frac{7}{30} c_1,  \
c_{12} = \frac{4}{5} c_1,  \\
&c_{13} = - \frac{13}{10} c_1,  \
c_{14} = - \frac{4}{5} c_1,  \
c_{15} = - \frac{4}{5} c_1
\end{split}
\]
and these are the coefficients of Getzler's equation in $\ocM_{2,2}$.

\subsection{Checking $\mathfrak{r}_2(E)=0$}
Again, one has to pick a basis and check all components vanish.

\[
\< \d^\mu \>_1 \< \d^j \d^x \d^y \> \< \d^i \d^\mu \d^\nu \> \< \d^\nu \d^\alpha \d^\alpha \>: \quad
 \frac{1}{20} c_2 + c_4 - c_5 - c_8 + \frac{1}{24} c_{12} = 0.
\]
\begin{equation*}
\< \d^j \d^x \d^y \> \< \d^i \d^\mu \d^\nu \d^\nu \> \< \d^\mu \d^\alpha \d^\alpha \> : \quad
 \frac{1}{1152} c_2 - c_3 + \frac{1}{24} c_4 - \frac{1}{24} c_8 = 0. 
\end{equation*} 
\begin{equation*}
\< \d^j \d^x \d^y \> \< \d^i \d^\mu \d^\nu \> \< \d^\mu \d^\nu \d^\alpha \d^\alpha \> : \quad
 \frac{1}{480} c_2 - 2 c_3 + \frac{1}{12} c_4 = 0. 
\end{equation*} 
\begin{equation*}  
\< \d^j \>_1 \< \d^x \d^y \d^\mu \> \< \d^i \d^\mu \d^\nu \> \< \d^\nu \d^\alpha \d^\alpha \> : \quad
 - \frac{1}{12} c_1 - \frac{1}{24} c_2 - c_4 + c_5 - c_8 - c_9 - c_{10} + \frac{1}{24} c_{12} = 0. 
\end{equation*} 
\begin{equation*}  \< \d^\mu \>_1 \< \d^j \d^\nu \d^\nu \> \< \d^i \d^\mu \d^\alpha \> \< \d^x \d^y \d^\alpha \> : \quad
 \frac{1}{12} c_1 + \frac{1}{24} c_2 - c_4 - c_5 + c_8 + \frac{1}{24} c_{12} + \frac{1}{24} c_{14} + \frac{1}{24} c_{15} = 0. 
\end{equation*} 
\begin{equation*}  \< \d^\mu \>_1 \< \d^j \d^y \d^\mu \> \< \d^i \d^\nu \d^\alpha \> \< \d^x \d^\nu \d^\alpha \> : \quad
 - \frac{1}{24} c_1 - c_5 + c_8 + c_{10} + \frac{1}{24} c_{15} = 0. 
\end{equation*} 
\begin{equation*}  \< \d^\mu \>_1 \< \d^j \d^x \d^\mu \> \< \d^i \d^\nu \d^\alpha \> \< \d^y \d^\nu \d^\alpha \> : \quad
 - \frac{1}{24} c_1 - c_5 + c_8 + c_9 + \frac{1}{24} c_{14} = 0. 
\end{equation*} 
\begin{equation*}  \< \d^j \d^x \d^\mu \d^\mu \> \< \d^i \d^\nu \d^\alpha \> \< \d^y \d^\nu \d^\alpha \> : \quad
 - \frac{1}{576} c_1 - 2 c_3 + \frac{1}{24} c_6 + \frac{1}{24} c_8 + \frac{1}{24} c_9 = 0. 
\end{equation*} 
\begin{equation*}  \< \d^j \d^y \d^\mu \d^\mu \> \< \d^i \d^\nu \d^\alpha \> \< \d^x \d^\nu \d^\alpha \> : \quad
 - \frac{1}{576} c_1 - 2 c_3 + \frac{1}{24} c_7 + \frac{1}{24} c_8 + \frac{1}{24} c_{10} = 0. 
\end{equation*}

The $7$ vectors
\[
\begin{split}
&\< \d^i \d^\mu \d^\mu \d^\nu \> \< \d^x \d^y \d^\nu \> \< \d^j \d^\alpha \d^\alpha \>,
\< \d^x \d^\mu \d^\nu \d^\nu \> \< \d^i \d^y \d^\mu \> \< \d^j \d^\alpha \d^\alpha \>, \\
&\< \d^y \d^\mu \d^\nu \d^\nu \> \< \d^i \d^x \d^\mu \> \< \d^j \d^\alpha \d^\alpha \>,
\< \d^x \d^y \d^\mu \d^\nu \> \< \d^i \d^\mu \d^\nu \> \< \d^j \d^\alpha \d^\alpha \> ,\\
&\< \d^i \d^x \d^y \d^\mu \> \< \d^\mu \d^\nu \d^\nu \> \< \d^j \d^\alpha \d^\alpha \> ,
\< \d^i \d^x \d^\mu \d^\nu \> \< \d^y \d^\mu \d^\nu \> \< \d^j \d^\alpha \d^\alpha \> , \\
&\< \d^i \d^y \d^\mu \d^\nu \> \< \d^x \d^\mu \d^\nu \> \< \d^j \d^\alpha \d^\alpha \>
\end{split}
\]
are related by WDVV equations.
$4$ of them are linear independent.
\[
 \< \d^i \d^\mu \d^\mu \d^\nu \> \< \d^x \d^y \d^\nu \> \< \d^j \d^\alpha \d^\alpha \>  : \quad
 \frac{1}{288} c_1 + \frac{1}{576} c_2 - 2 c_3 + \frac{1}{12} c_8 + \frac{1}{24} c_9 + \frac{1}{24} c_{10} = 0. 
\]
\[
 \< \d^x \d^\mu \d^\nu \d^\nu \> \< \d^i \d^y \d^\mu \> \< \d^j \d^\alpha \d^\alpha \>  : \quad
 \frac{1}{24} c_6 - \frac{1}{24} c_8 - \frac{1}{24} c_9 = 0. 
\]
\[
 \< \d^y \d^\mu \d^\nu \d^\nu \> \< \d^i \d^x \d^\mu \> \< \d^j \d^\alpha \d^\alpha \>  : \quad
 \frac{1}{24} c_7 - \frac{1}{24} c_8 - \frac{1}{24} c_{10} = 0. 
\]
\[
 \< \d^x \d^y \d^\mu \d^\nu \> \< \d^i \d^\mu \d^\nu \> \< \d^j \d^\alpha \d^\alpha \>  : \quad
 \frac{1}{288} c_1 - 4 c_3 + \frac{1}{6} c_8 + \frac{1}{24} c_9 + \frac{1}{24} c_{10} + \frac{1}{24} c_{11} = 0. 
\]

All of the remaining strata are related by WDVV equations. 
Only $3$ of them are linearly independent.
\[
 \< \d^x \d^j \d^\mu \> \< \d^\mu \d^y \d^i \d^\nu \> \< \d^\nu \d^\alpha \d^\alpha \>  : \quad
 - 4 c_3 - \frac{5}{6} c_4 + c_8 + \frac{1}{6} c_9 + \frac{1}{6} c_{10} - \frac{1}{12} c_{11} = 0. 
\]
\[
 \< \d^x \d^j \d^\alpha \d^\mu \> \< \d^\mu \d^i \d^\nu \> \< \d^\nu \d^y \d^\alpha \>  : \quad
 - \frac{5}{6} c_4 - \frac{5}{6} c_9 = 0. 
\]
\[
 \< \d^y \d^j \d^\alpha \d^\mu \> \< \d^\mu \d^i \d^\nu \> \< \d^\nu \d^x \d^\alpha \>  : \quad
 - \frac{5}{6} c_4 - \frac{5}{6} c_{10} = 0. 
\]

Therefore, $\mathfrak{r}_2 (E)=0$.

\subsection{Calculating $\mathfrak{r}_3 (E)$}
Since $l=3$ case is new, the calculation is presented.

\begin{eqnarray*}
& & \< \d^x_1 \d^y_1 \>_2 
\mapsto \\ & & \hspace{.5cm} 
- \frac{1}{24}  \< \d^y \d^\nu \d^j \> \< \d^x \d^\alpha \d^\mu \> \< \d^\alpha \d^\nu \d^i \> \< \d^\mu \d^\beta \d^\beta \>
- \frac{1}{24}  \< \d^y \d^\nu \d^\mu \> \< \d^x \d^\alpha \d^j \> \< \d^\alpha \d^\nu \d^i \> \< \d^\mu \d^\beta \d^\beta \>
\\ & & \hspace{.5cm} 
- \frac{1}{24}  \< \d^x \d^\alpha \d^j \> \< \d^\alpha \d^i \d^\mu \> \< \d^y \d^\beta \d^\nu \> \< \d^\beta \d^\mu \d^\nu \> 
- \frac{1}{24}  \< \d^y \d^\alpha \d^j \> \< \d^\alpha \d^i \d^\mu \> \< \d^x \d^\beta \d^\nu \> \< \d^\beta \d^\mu \d^\nu \> 
\\ & & \hspace{.5cm} 
- \frac{1}{24}  \< \d^i \d^\mu \d^j \> \< \d^y \d^\beta \d^\nu \> \< \d^x \d^\alpha \d^\nu \> \< \d^\alpha \d^\beta \d^\mu \> 
- \frac{1}{24}  \< \d^i \d^\mu \d^j \> \< \d^y \d^\beta \d^\nu \> \< \d^x \d^\alpha \d^\nu \> \< \d^\alpha \d^\beta \d^\mu \> 
\\ & & \hspace{.5cm} 
+ \frac{1}{48}  \< \d^\alpha \d^y \d^\beta \> \< \d^\beta \d^x \d^\mu \> \< \d^\mu \d^j \d^\nu \> \< \d^\nu \d^i \d^\alpha \> 
+ \frac{1}{48}  \< \d^x \d^\mu \d^\alpha \> \< \d^\nu \d^y \d^\beta \> \< \d^\beta \d^\mu \d^j \> \< \d^\nu \d^i \d^\alpha \> 
\\ & & \hspace{.5cm} 
+ \frac{1}{48}  \< \d^x \d^\mu \d^\alpha \> \< \d^\mu \d^j \d^\nu \> \< \d^\alpha \d^y \d^\beta \> \< \d^\beta \d^\nu \d^i \> 
+ \frac{1}{48}  \< \d^\alpha \d^y \d^\beta \> \< \d^\beta \d^j \d^\nu \> \< \d^x \d^\mu \d^\alpha \> \< \d^\mu \d^\nu \d^i \> 
\\ & & \hspace{.5cm}
+ \frac{1}{48}  \< \d^j \d^\nu \d^\alpha \> \< \d^\alpha \d^y \d^\beta \> \< \d^\beta \d^x \d^\mu \> \< \d^\mu \d^\nu \d^i \> 
+ \frac{1}{48}  \< \d^j \d^\nu \d^\alpha \> \< \d^x \d^\mu \d^\alpha \> \< \d^i \d^y \d^\beta \> \< \d^\beta \d^\mu \d^\nu \> 
\\ & & \hspace{.5cm}
+ \frac{1}{576}  \< \d^y \d^\mu \d^\alpha \> \< \d^\mu \d^x \d^\nu \> \< \d^\nu \d^i \d^\alpha \> \< \d^j \d^\beta \d^\beta \> 
+ \frac{1}{576}  \< \d^x \d^\nu \d^\alpha \> \< \d^y \d^\mu \d^\alpha \> \< \d^\mu \d^\nu \d^i \> \< \d^j \d^\beta \d^\beta \> 
\\ & & \hspace{.5cm}
+ \frac{1}{576}  \< \d^x \d^\mu \d^\nu \> \< \d^\mu \d^i \d^\nu \> \< \d^y \d^\alpha \d^\beta \> \< \d^\alpha \d^j \d^\beta \> \\
\end{eqnarray*}
\begin{eqnarray*}
& & \< \d^\mu_1 \>_2 \< \d^x \d^y \d^\mu \> 
\mapsto \\ & & \hspace{.5cm} 
- \frac{1}{5760}  \< \d^i \d^\alpha \d^\mu \> \< \d^\alpha \d^\beta \d^\beta \> \< \d^\mu \d^\nu \d^\nu \> \< \d^x \d^y \d^j \> 
+ \frac{1}{960}  \< \d^i \d^\alpha \d^\nu \> \< \d^\alpha \d^\beta \d^\nu \> \< \d^\beta \d^\mu \d^\mu \> \< \d^x \d^y \d^j \> 
\\ & & \hspace{.5cm} 
- \frac{1}{24}  \< \d^\mu \d^j \d^\beta \> \< \d^\beta \d^i \d^\nu \> \< \d^\nu \d^\alpha \d^\alpha \> \< \d^x \d^y \d^\mu \> 
- \frac{1}{24}  \< \d^i \d^\nu \d^j \> \< \d^\alpha \d^\mu \d^\beta \> \< \d^\beta \d^\nu \d^\alpha \> \< \d^x \d^y \d^\mu \> 
\\ & & \hspace{.5cm} 
+ \frac{1}{48}  \< \d^\mu \d^\beta \d^\alpha \> \< \d^\beta \d^j \d^\nu \> \< \d^\nu \d^i \d^\alpha \> \< \d^x \d^y \d^\mu \> 
+ \frac{1}{48}  \< \d^j \d^\nu \d^\alpha \> \< \d^\mu \d^\beta \d^\alpha \> \< \d^\beta \d^\nu \d^i \> \< \d^x \d^y \d^\mu \> 
\\ & & \hspace{.5cm} 
+ \frac{1}{576}  \< \d^\mu \d^\nu \d^\alpha \> \< \d^\nu \d^i \d^\alpha \> \< \d^j \d^\beta \d^\beta \> \< \d^x \d^y \d^\mu \> \\
\end{eqnarray*}
\begin{eqnarray*}
& & \< \d^x \d^y \d^\mu \d^\mu \d^\nu \d^\nu \> 
\mapsto \\ & & \hspace{.5cm} 
+ 4  \< \d^i \d^\nu \d^\beta \> \< \d^\beta \d^\mu \d^\nu \> \< \d^\mu \d^\alpha \d^j \> \< \d^\alpha \d^x \d^y \>
- 2  \< \d^i \d^\alpha \d^\mu \> \< \d^\alpha \d^\beta \d^\mu \> \< \d^\beta \d^x \d^y \> \< \d^j \d^\nu \d^\nu \> 
\\ & & \hspace{.5cm} 
- 4  \< \d^i \d^\alpha \d^\nu \> \< \d^\alpha \d^\beta \d^\mu \> \< \d^\beta \d^x \d^y \> \< \d^j \d^\mu \d^\nu \> 
- 4  \< \d^i \d^\alpha \d^\nu \> \< \d^\alpha \d^\beta \d^\mu \> \< \d^\beta \d^x \d^\mu \> \< \d^j \d^y \d^\nu \> 
\\ & & \hspace{.5cm} 
- 4  \< \d^i \d^\alpha \d^\nu \> \< \d^\alpha \d^\beta \d^\mu \> \< \d^\beta \d^y \d^\mu \> \< \d^j \d^x \d^\nu \> 
-  \< \d^i \d^\alpha \d^\nu \> \< \d^\alpha \d^\beta \d^\nu \> \< \d^\beta \d^\mu \d^\mu \> \< \d^j \d^x \d^y \> 
\\ & & \hspace{.5cm} 
+ 4  \< \d^i \d^\mu \d^\alpha \> \< \d^\alpha \d^x \d^y \> \< \d^j \d^\nu \d^\beta \> \< \d^\beta \d^\mu \d^\nu \> 
+ 2  \< \d^i \d^\mu \d^\alpha \> \< \d^\alpha \d^x \d^\mu \> \< \d^j \d^\nu \d^\beta \> \< \d^\beta \d^y \d^\nu \> 
\\ & & \hspace{.5cm} 
+ 4  \< \d^i \d^\nu \d^\alpha \> \< \d^\alpha \d^x \d^\mu \> \< \d^j \d^\nu \d^\beta \> \< \d^\beta \d^y \d^\mu \> \\
\end{eqnarray*}
\begin{eqnarray*}
& & \< \d^\mu \d^\mu \d^\nu \>_1 \< \d^x \d^y \d^\nu \> 
\mapsto \\ & & \hspace{.5cm} 
+ \frac{1}{12}  \< \d^i \d^\nu \d^\beta \> \< \d^\beta \d^\mu \d^j \> \< \d^\mu \d^\alpha \d^\alpha \> \< \d^x \d^y \d^\nu \> 
+ \frac{1}{24}  \< \d^i \d^\mu \d^\beta \> \< \d^\beta \d^\nu \d^\mu \> \< \d^\nu \d^\alpha \d^\alpha \> \< \d^x \d^y \d^j \> 
\\ & & \hspace{.5cm} 
-  \< \d^i \d^\alpha \d^\nu \> \< \d^\alpha \d^\beta \d^\mu \> \< \d^\beta \d^j \d^\mu \> \< \d^x \d^y \d^\nu \> 
+ \frac{1}{2}  \< \d^i \d^\beta \d^\mu \> \< \d^j \d^\alpha \d^\nu \> \< \d^\alpha \d^\beta \d^\mu \> \< \d^x \d^y \d^\nu \> 
\\ & & \hspace{.5cm} 
+ \frac{1}{2}  \< \d^i \d^\beta \d^\nu \> \< \d^j \d^\alpha \d^\mu \> \< \d^\alpha \d^\beta \d^\mu \> \< \d^x \d^y \d^\nu \> 
- \frac{1}{12}  \< \d^i \d^\alpha \d^\mu \> \< \d^\alpha \d^\beta \d^\beta \> \< \d^j \d^\mu \d^\nu \> \< \d^x \d^y \d^\nu \> 
\\ & & \hspace{.5cm} 
- \frac{1}{24}  \< \d^i \d^\alpha \d^\nu \> \< \d^\alpha \d^\beta \d^\beta \> \< \d^j \d^\mu \d^\mu \> \< \d^x \d^y \d^\nu \> 
+ \frac{1}{24}  \< \d^i \d^\alpha \d^\alpha \> \< \d^j \d^\beta \d^\nu \> \< \d^\beta \d^\mu \d^\mu \> \< \d^x \d^y \d^\nu \> \\
\end{eqnarray*}
\begin{eqnarray*}
& & \< \d^\mu \d^\mu \d^\nu \> \< \d^x \d^y \d^\nu \>_1 
\mapsto \\ & & \hspace{.5cm} 
+ \frac{1}{24}  \< \d^\mu \d^\mu \d^j \> \< \d^i \d^y \d^\beta \> \< \d^\beta \d^\nu \d^x \> \< \d^\nu \d^\alpha \d^\alpha \>
-  \< \d^\mu \d^\mu \d^\nu \> \< \d^i \d^\alpha \d^\nu \> \< \d^\alpha \d^\beta \d^y \> \< \d^\beta \d^j \d^x \>
\\ & & \hspace{.5cm} 
+ \frac{1}{2}  \< \d^\mu \d^\mu \d^\nu \> \< \d^i \d^\beta \d^y \> \< \d^j \d^\alpha \d^\nu \> \< \d^\alpha \d^\beta \d^x \>
+ \frac{1}{2}  \< \d^\mu \d^\mu \d^\nu \> \< \d^i \d^\beta \d^\nu \> \< \d^j \d^\alpha \d^y \> \< \d^\alpha \d^\beta \d^x \>
\\ & & \hspace{.5cm} 
- \frac{1}{24}  \< \d^\mu \d^\mu \d^\nu \> \< \d^i \d^\alpha \d^x \> \< \d^\alpha \d^\beta \d^\beta \> \< \d^j \d^y \d^\nu \> 
- \frac{1}{24}  \< \d^\mu \d^\mu \d^\nu \> \< \d^i \d^\alpha \d^y \> \< \d^\alpha \d^\beta \d^\beta \> \< \d^j \d^x \d^\nu \> 
\\ & & \hspace{.5cm} 
- \frac{1}{24}  \< \d^\mu \d^\mu \d^\nu \> \< \d^i \d^\alpha \d^\nu \> \< \d^\alpha \d^\beta \d^\beta \> \< \d^j \d^x \d^y \> 
+ \frac{1}{24}  \< \d^\mu \d^\mu \d^\nu \> \< \d^i \d^\alpha \d^\alpha \> \< \d^j \d^\beta \d^\nu \> \< \d^\beta \d^x \d^y \> \\
\end{eqnarray*}
\begin{eqnarray*}
& & \< \d^x \d^\mu \d^\nu \>_1 \< \d^y \d^\mu \d^\nu \> 
\mapsto \\ & & \hspace{.5cm} 
+ \frac{1}{12}  \< \d^i \d^\nu \d^\beta \> \< \d^\beta \d^\mu \d^x \> \< \d^\mu \d^\alpha \d^\alpha \> \< \d^y \d^j \d^\nu \> 
-  \< \d^i \d^\alpha \d^\nu \> \< \d^\alpha \d^\beta \d^\mu \> \< \d^\beta \d^j \d^x \> \< \d^y \d^\mu \d^\nu \> 
\\ & & \hspace{.5cm} 
+ \frac{1}{2}  \< \d^i \d^\beta \d^\mu \> \< \d^j \d^\alpha \d^\nu \> \< \d^\alpha \d^\beta \d^x \> \< \d^y \d^\mu \d^\nu \> 
+ \frac{1}{2}  \< \d^i \d^\beta \d^\nu \> \< \d^j \d^\alpha \d^\mu \> \< \d^\alpha \d^\beta \d^x \> \< \d^y \d^\mu \d^\nu \> 
\\ & & \hspace{.5cm} 
- \frac{1}{24}  \< \d^i \d^\alpha \d^x \> \< \d^\alpha \d^\beta \d^\beta \> \< \d^j \d^\mu \d^\nu \> \< \d^y \d^\mu \d^\nu \> 
- \frac{1}{12}  \< \d^i \d^\alpha \d^\mu \> \< \d^\alpha \d^\beta \d^\beta \> \< \d^j \d^x \d^\nu \> \< \d^y \d^\mu \d^\nu \> 
\\ & & \hspace{.5cm} 
+ \frac{1}{24}  \< \d^i \d^\alpha \d^\alpha \> \< \d^j \d^\nu \d^\beta \> \< \d^\beta \d^x \d^\mu \> \< \d^y \d^\mu \d^\nu \> \\
\end{eqnarray*}
\begin{eqnarray*}
& & \< \d^y \d^\mu \d^\nu \>_1 \< \d^x \d^\mu \d^\nu \> 
\mapsto \\ & & \hspace{.5cm} 
+ \frac{1}{12}  \< \d^i \d^\nu \d^\beta \> \< \d^\beta \d^\mu \d^y \> \< \d^\mu \d^\alpha \d^\alpha \> \< \d^x \d^j \d^\nu \> 
-  \< \d^i \d^\alpha \d^\nu \> \< \d^\alpha \d^\beta \d^\mu \> \< \d^\beta \d^j \d^y \> \< \d^x \d^\mu \d^\nu \> 
\\ & & \hspace{.5cm} 
+ \frac{1}{2}  \< \d^i \d^\beta \d^\mu \> \< \d^j \d^\alpha \d^\nu \> \< \d^\alpha \d^\beta \d^y \> \< \d^x \d^\mu \d^\nu \> 
+ \frac{1}{2}  \< \d^i \d^\beta \d^\nu \> \< \d^j \d^\alpha \d^\mu \> \< \d^\alpha \d^\beta \d^y \> \< \d^x \d^\mu \d^\nu \> 
\\ & & \hspace{.5cm} 
- \frac{1}{24}  \< \d^i \d^\alpha \d^y \> \< \d^\alpha \d^\beta \d^\beta \> \< \d^j \d^\mu \d^\nu \> \< \d^x \d^\mu \d^\nu \> 
- \frac{1}{12}  \< \d^i \d^\alpha \d^\mu \> \< \d^\alpha \d^\beta \d^\beta \> \< \d^j \d^y \d^\nu \> \< \d^x \d^\mu \d^\nu \> 
\\ & & \hspace{.5cm} 
+ \frac{1}{24}  \< \d^i \d^\alpha \d^\alpha \> \< \d^j \d^\nu \d^\beta \> \< \d^\beta \d^y \d^\mu \> \< \d^x \d^\mu \d^\nu \> \\
\end{eqnarray*}
The other graphs all have $\mathfrak{r}_3(\Gamma)=0$.

\subsection{Checking $\mathfrak{r}_3 (E)=0$}

Now the $\mathfrak{r}_3(E)$ has only a few independent strata.

\[
\< \d^i \d^\mu \d^\mu \> : \quad
 \frac{1}{288} c_1 + \frac{1}{576} c_2 - 2 c_3 + \frac{1}{12} c_8 + \frac{1}{24} c_9 + \frac{1}{24} c_{10} = 0. 
\]
\[
\< \d^i \d^x \d^y \>: \quad 
 \frac{1}{1152} c_2 - c_3 + \frac{1}{24} c_4 - \frac{1}{24} c_8 = 0. 
\]
\[
\< \d^i \d^x \d^\mu \> \< \d^\mu \d^\nu \d^\nu \> \< \d^j \d^y \d^\alpha \> \< \d^\alpha \d^\beta \d^\beta \>: \quad
 \frac{1}{576} c_1 + 2 c_3 - \frac{1}{12} c_8 - \frac{1}{24} c_9 - \frac{1}{24} c_{10} = 0. 
\]
The remaining strata, which are all equivalent as codimension $4$ classes in
$\ocM_{1,4}$, have coefficient
$$ - \frac{1}{8} c_1 - \frac{1}{24} c_2 = 0. $$

\section{The Belorousski-Pandharipande equation in $\ocM_{2,3}$}

\subsection{Strata in $\ocM_{2,3}$ of codimension $2$}
\emph{Throughout this section, we will always assume that the three 
external labelings $x,y,z$ are symmetrized.}

Using Mumford--Getzler's and Getzler's equations for genus $2$, 
all genus $2$ terms
with more than one descendent can be rewritten in terms of the others.
Also, using TRR's, all genus $0$ or $1$ terms with a descendent can be
rewritten in terms of the others.
We are therefore left with $21$ strata which could be independent.
\footnote{We
rewrote the strata here in the same order as in the Belorousski-Pandharipande
equation, with the extra strata
$\< \d^x \d^y \d^z \d^\mu \d^\mu \d^\nu \d^\nu \>$, which does not appear in
the equation, at the end.} 

A general linear combination is
\begin{eqnarray*}
E = \sum_{S_3 (x,y,z)} 
& &  c_1 \< \d^x \d^\mu \d^\nu \> \< \d^\mu \d^y \d^z \> \< \d^\nu \>_2
+ c_2 \< \d^x \d^y \d^z \d^\mu \> \< \d^\mu_1 \>_2
+ c_3 \< \d^x_1 \d^\mu \>_2 \< \d^\mu \d^y \d^z \>  \\
& & + c_4 \< \d^x \d^\mu_1 \>_2 \< \d^\mu \d^y \d^z \>
+ c_5 \< \d^x \d^y \d^z \d^\mu \d^\nu \> \< \d^\mu \>_1 \< \d^\nu \>_1 \\
& & + c_6 \< \d^x \d^y \d^\mu \d^\nu \> \< \d^\mu \>_1 \< \d^\nu \d^z \>_1
+ c_7 \< \d^x \d^\mu \d^\nu \> \< \d^\mu \>_1 \< \d^\nu \d^y \d^z \>_1 \\
& & + c_8 \< \d^x \d^\mu \>_1 \< \d^\mu \d^y \d^\nu \> \< \d^\nu \d^z \>_1
+ c_9 \< \d^x \d^y \d^z \d^\mu \> \< \d^\mu \d^\nu \>_1 \< \d^\nu \>_1 \\
& & + c_{10} \< \d^x \d^\mu \>_1 \< \d^\mu \d^\nu \>_1 \< \d^\nu \d^y \d^z \>
+ c_{11} \< \d^\mu \>_1 \< \d^\mu \d^x \d^\nu \>_1 \< \d^\nu \d^y \d^z \> \\
& & + c_{12} \< \d^x \d^y \d^z \d^\mu \d^\mu \d^\nu \> \< \d^\nu \>_1
+ c_{13} \< \d^x \d^y \d^\mu \d^\mu \d^\nu \> \< \d^\nu \d^z \>_1 \\
& & + c_{14} \< \d^x \d^\mu \d^\mu \d^\nu \> \< \d^\nu \d^y \d^z \>_1
+ c_{15} \< \d^x \d^y \d^z \d^\mu \>_1 \< \d^\mu \d^\nu \d^\nu \> \\
& & + c_{16} \< \d^x \d^y \d^z \d^\mu \d^\nu \> \< \d^\mu \d^\nu \>_1
+ c_{17} \< \d^x \d^y \d^\mu \d^\nu \> \< \d^\mu \d^\nu \d^z \>_1 \\
& & + c_{18} \< \d^x \d^\mu \d^\nu \> \< \d^\mu \d^\nu \d^y \d^z \>_1
+ c_{19} \< \d^x \d^y \d^z \d^\mu \> \< \d^\mu \d^\nu \d^\nu \>_1 \\
& & + c_{20} \< \d^x \d^y \d^\mu \> \< \d^\mu \d^\nu \d^\nu \d^z \>_1
+ c_{21} \< \d^x \d^y \d^z \d^\mu \d^\mu \d^\nu \d^\nu \> = 0,
\end{eqnarray*}

\subsection{Setting $\mathfrak{r}_1(E)=0$}
Again a basis is chosen for the output graphs of $\mathfrak{r}_1(E)$, and 
the coefficients are set to zero.

\begin{equation}\label{eq231}
\< \d^x \d^\mu \d^j \> \< \d^\mu \d^y \d^z \> \< \d^i_1 \>_2 : \quad
c_1 + c_2 = 0.
\end{equation}
\begin{equation}\label{eq232}
\< \d^x \d^y \d^i \> \< \d^z \d^\mu \d^j \> \< \d^\mu_1 \>_2: \quad
- 3 c_2 + 3 c_3 + c_4 = 0.
\end{equation}
\begin{equation}\label{eq233}
\< \d^x \d^y \d^i \>_1 \< \d^z \d^j \d^\mu \> \< \d^\mu \d^\nu \d^\nu \>: \quad
c_{14} - 3 c_{15} - c_{18} = 0.
\end{equation}
\begin{equation}\label{eq234}
\< \d^x \d^y \d^\mu \>_1 \< \d^z \d^i \d^\mu \> \< \d^j \d^\nu \d^\nu \>: \quad
\frac{1}{24} c_7 - c_{14} + c_{15} + c_{15} + c_{15} = 0.
\end{equation}
\begin{equation}\label{eq235}
\< \d^x \d^i \d^\mu \>_1 \< \d^\mu \d^\nu \d^\nu \> \< \d^y \d^z \d^j \>: \quad
- \frac{1}{80} c_3 - \frac{1}{240} c_4 - 3 c_{15} = 0.
\end{equation}
\begin{equation}\label{eq236}
\< \d^x \d^i \d^\mu \>_1 \< \d^y \d^z \d^\mu \> \< \d^j \d^\nu \d^\nu \>: \quad
- c_{20} + \frac{1}{24} c_{11} = 0.
\end{equation}
\begin{equation}\label{eq237}
\< \d^x \d^\mu \d^\nu \>_1 \< \d^i \d^\mu \d^\nu \> \< \d^y \d^z \d^j \>: \quad
\frac{1}{30} c_3 + \frac{1}{10} c_4 - c_{17} + c_{20} + c_{20} = 0.
\end{equation}
\begin{equation}\label{eq238}
\< \d^i \d^\mu \d^\nu \>_1 \< \d^x \d^\mu \d^\nu \> \< \d^y \d^z \d^j \>: \quad
\frac{1}{30} c_3 - c_{18} = 0.
\end{equation}
\begin{equation}\label{eq239}
\< \d^i \d^\mu \d^\mu \>_1 \< \d^x \d^y \d^\nu \> \< \d^\nu \d^z \d^j \>: \quad
c_{19} - c_{20} = 0.
\end{equation}
\begin{equation}\label{eq2310}
\< \d^\mu \d^\nu \d^\nu \>_1 \< \d^x \d^i \d^\mu \> \< \d^y \d^z \d^j \>: \quad
- \frac{1}{30} c_3 - 2 c_{19} - c_{19} + c_{20} = 0.
\end{equation}
\begin{equation}\label{eq2311}
\< \d^x \d^i \>_1 \< \d^\mu \>_1 \< \d^y \d^z \d^\nu \> \< \d^\nu \d^j \d^\mu \>: \quad
- c_4 + c_6 - 2 c_7 + c_{10} - c_{11} = 0.
\end{equation}
\begin{equation}\label{eq2312}
\< \d^x \d^\mu \>_1 \< \d^\nu \>_1 \< \d^i \d^\mu \d^\nu \> \< \d^y \d^z \d^j \>: \quad 
\frac{4}{5} c_3 + \frac{7}{5} c_4 - c_6 + c_{10} + c_{11} = 0.
\end{equation}
\begin{equation}\label{eq2313}
\< \d^x \d^\mu \>_1 \< \d^\mu \d^y \d^i \> \< \d^\nu \>_1 \< \d^\nu \d^z \d^j \>: \quad 
- 2 c_6 + c_7 + c_7 + 2 c_8 = 0.
\end{equation}
\begin{equation}\label{eq2314}
\< \d^i \d^\mu \>_1 \< \d^\nu \>_1 \< \d^x \d^\mu \d^\nu \> \< \d^y \d^z \d^j \>: \quad 
\frac{4}{5} c_3 - c_7 = 0.
\end{equation}
\begin{equation}\label{eq2315}
\< \d^i \d^\mu \>_1 \< \d^\mu \d^x \d^y \> \< \d^\nu \>_1 \< \d^\nu \d^z \d^j \>: \quad 
- c_3 + c_{10} - c_{11} = 0.
\end{equation}
\begin{equation}\label{eq2316}
\< \d^i \d^\mu \>_1 \< \d^\mu \>_1 \< \d^x \d^y \d^\nu \> \< \d^\nu \d^z \d^j \>: \quad 
c_9 - c_{11} = 0.
\end{equation}
\begin{equation}\label{eq2317}
\< \d^\mu \d^\nu \>_1 \< \d^\mu \>_1 \< \d^x \d^i \d^\nu \> \< \d^y \d^z \d^j \>: \quad 
- \frac{4}{5} c_3 - 2 c_9 - c_9 + c_{11} = 0.
\end{equation}
\begin{equation}\label{eq2318}
\< \d^\mu \d^\nu \>_1 \< \d^i \d^\mu \d^\nu \> \< \d^x \d^y \d^z \d^j \>: \quad ,
\frac{1}{10} c_2 - c_{16} + c_{19} + c_{19} = 0.
\end{equation}
\begin{equation}\label{eq2319}
\< \d^i \d^\mu \>_1 \< \d^\mu \d^\nu \d^\nu \> \< \d^x \d^y \d^z \d^j \>: \quad ,
- \frac{1}{240} c_2 - c_{15} = 0.
\end{equation}
\begin{equation}\label{eq2320}
\< \d^x \d^\mu \>_1 \< \d^i \d^\mu \d^\nu \d^\nu \> \< \d^y \d^z \d^j \>: \quad
\frac{1}{48}c_3 + \frac{13}{240}c_4 + \frac{1}{24} c_{10} - c_{13} + c_{20}=0.
\end{equation}
\begin{equation}\label{eq2321}
\< \d^\mu \>_1 \< \d^\nu \>_1 \< \d^i \d^\mu \d^\nu \> \< \d^x \d^y \d^z \d^j \>: \quad 
\frac{7}{10} c_2 - c_5 + c_9 = 0.
\end{equation}
\begin{equation}\label{eq2322}
\< \d^\mu \>_1 \< \d^i \d^\mu \d^\nu \d^\nu \> \< \d^x \d^y \d^z \d^j \>: \quad 
\frac{13}{240} c_2 + \frac{1}{24} c_9 - c_{12} + c_{19} = 0.
\end{equation}
\begin{equation}\label{eq2323}
\< \d^x \d^y \d^i \> \< \d^j \d^z \d^\mu \d^\mu \d^\nu \d^\nu \>: \quad 
\frac{1}{576} c_3 + \frac{1}{960} c_4 + \frac{1}{24} c_{20} + 3 c_{21} = 0.
\end{equation}
\begin{equation}\label{eq2324}
\< \d^\mu \d^\mu \d^i \> \< \d^j \d^x \d^y \d^z \d^\nu \d^\nu \>: \quad 
\frac{1}{24} c_{12} + \frac{1}{24} c_{15} + 2 c_{21} = 0.
\end{equation}
\begin{equation}\label{eq2325}
\< \d^x \d^y \d^z \d^i \> \< \d^j \d^\mu \d^\mu \d^\nu \d^\nu \>: \quad 
\frac{1}{960} c_2 + \frac{1}{24} c_{19} + c_{21} = 0.
\end{equation}
\begin{equation}\label{eq2326}
\< \d^x \d^\mu \d^\mu \d^i \> \< \d^j \d^y \d^z \d^\nu \d^\nu \>: \quad 
\frac{1}{24} c_{13} + \frac{1}{24} c_{14} + 6 c_{21} = 0.
\end{equation}

The coefficients of the Belorousski-Pandharipande equation are the only
solution of the equations (\ref{eq231})-(\ref{eq2326}).

\subsection{Checking $\mathfrak{r}_2(E)=0$}

Since the whole output graphs are far too numerous, we shall present the
the coefficients of the following four (disconnected) graphs are zero:
\[
 \begin{split}
 \< \d^i \d^x \>_1 \< \d^j \d^y \d^\mu \> \< \d^\mu \d^z \d^\nu \> \< \d^\nu \d^\alpha \d^\alpha \>, \quad 
\< \d^i \d^\mu \>_1 \< \d^\mu \d^x \d^y \> \< \d^j \d^z \d^\nu \> \< \d^\nu \d^\alpha \d^\alpha \> \\
 \< \d^\mu \d^\nu \>_1 \< \d^i \d^\mu \d^\nu \> \< \d^j \d^x \d^\alpha \> \< \d^\alpha \d^y \d^z \>, \quad 
\< \d^x \d^\mu \>_1 \< \d^\mu \d^i \d^y \> \< \d^j \d^z \d^\nu \> \< \d^\nu \d^\alpha \d^\alpha \> 
 \end{split}
\]

In the following, the graphs $\Gamma$ appearing in BP equation such that
$\mathfrak{r}_2(\Gamma)$ contain any of the above four graphs will be listed.

\[
\< \d^x \d^\mu \d^\nu \> \< \d^\mu \d^y \d^z \> \< \d^\nu \>_2 
\mapsto 
\frac{1}{10} \< \d^\mu \d^\nu \>_1 \< \d^i \d^\mu \d^\nu \> \< \d^j \d^x \d^\alpha \> \< \d^\alpha \d^y \d^z \> 
+ \ldots .
\]
\[
\< \d^x_1 \d^\mu \>_2 \< \d^\mu \d^y \d^z \> 
\mapsto
\frac{1}{24} \< \d^i \d^\mu \>_1 \< \d^\mu \d^x \d^y \> \< \d^j \d^z \d^\nu \> \< \d^\nu \d^\alpha \d^\alpha \> 
+ \ldots 
\]
\[
\< \d^x \d^\mu_1 \>_2 \< \d^\mu \d^y \d^z \> 
\mapsto
\frac{1}{24} \< \d^i \d^x \>_1 \< \d^j \d^y \d^\mu \> \< \d^\mu \d^z \d^\nu \> \< \d^\nu \d^\alpha \d^\alpha \> 
+ \dots
\] 
\[
\< \d^x \d^\mu \>_1 \< \d^\mu \d^y \d^\nu \> \< \d^\nu \d^z \>_1 
\mapsto
- \frac{1}{12} \< \d^x \d^\mu \>_1 \< \d^\mu \d^i \d^y \> \< \d^j \d^z \d^\nu \> \< \d^\nu \d^\alpha \d^\alpha \> 
+ \dots 
\]
\begin{eqnarray*} & & 
\< \d^x \d^\mu \>_1 \< \d^\mu \d^\nu \>_1 \< \d^\nu \d^y \d^z \> 
\mapsto \\ & & \hspace{.5cm}  
- \frac{1}{24} \< \d^i \d^x \>_1 \< \d^j \d^y \d^\mu \> \< \d^\mu \d^z \d^\nu \> \< \d^\nu \d^\alpha \d^\alpha \> 
- \frac{1}{24} \< \d^i \d^\mu \>_1 \< \d^\mu \d^x \d^y \> \< \d^j \d^z \d^\nu \> \< \d^\nu \d^\alpha \d^\alpha \> 
+ \dots 
\end{eqnarray*}
\begin{eqnarray*} & & 
\< \d^x \d^y \d^\mu \d^\mu \d^\nu \> \< \d^\nu \d^z \>_1 
\mapsto \\ & & \hspace{.5cm}  
- \< \d^i \d^x \>_1 \< \d^j \d^y \d^\mu \> \< \d^\mu \d^z \d^\nu \> \< \d^\nu \d^\alpha \d^\alpha \> 
+2 \< \d^x \d^\mu \>_1 \< \d^\mu \d^i \d^y \> \< \d^j \d^z \d^\nu \> \< \d^\nu \d^\alpha \d^\alpha \> 
+ \dots 
\end{eqnarray*}
\begin{eqnarray*} & & 
\< \d^x \d^y \d^z \d^\mu \>_1 \< \d^\mu \d^\nu \d^\nu \> 
\mapsto \\ & & \hspace{.5cm}  
-6 \< \d^x \d^\mu \>_1 \< \d^\mu \d^i \d^y \> \< \d^j \d^z \d^\nu \> \< \d^\nu \d^\alpha \d^\alpha \> 
+3 \< \d^i \d^x \>_1 \< \d^j \d^y \d^\mu \> \< \d^\mu \d^z \d^\nu \> \< \d^\nu \d^\alpha \d^\alpha \> 
+ \dots 
\end{eqnarray*}
\[
\< \d^x \d^y \d^z \d^\mu \d^\nu \> \< \d^\mu \d^\nu \>_1 
\mapsto
\< \d^\mu \d^\nu \>_1 \< \d^i \d^\mu \d^\nu \> \< \d^j \d^x \d^\alpha \> \< \d^\alpha \d^y \d^z \> 
+ \dots
\] 
\begin{eqnarray*} & & 
\< \d^x \d^\mu \d^\nu \> \< \d^\mu \d^\nu \d^y \d^z \>_1 \mapsto \\ & & \hspace{.5cm} 
2 \< \d^i \d^x \>_1 \< \d^j \d^y \d^\mu \> \< \d^\mu \d^z \d^\nu \> \< \d^\nu \d^\alpha \d^\alpha \> 
-2 \< \d^x \d^\mu \>_1 \< \d^\mu \d^i \d^y \> \< \d^j \d^z \d^\nu \> \< \d^\nu \d^\alpha \d^\alpha \> 
+ \dots 
\end{eqnarray*}
\begin{eqnarray*} & & 
\< \d^x \d^y \d^\mu \> \< \d^\mu \d^\nu \d^\nu \d^z \>_1 \mapsto \\ & & \hspace{.5cm} 
-2 \< \d^\mu \d^\nu \>_1 \< \d^i \d^\mu \d^\nu \> \< \d^j \d^x \d^\alpha \> \< \d^\alpha \d^y \d^z \> 
+ \< \d^i \d^x \>_1 \< \d^j \d^y \d^\mu \> \< \d^\mu \d^z \d^\nu \> \< \d^\nu \d^\alpha \d^\alpha \> 
\\ & & \hspace{.5cm} 
+ \< \d^i \d^\mu \>_1 \< \d^\mu \d^x \d^y \> \< \d^j \d^z \d^\nu \> \< \d^\nu \d^\alpha \d^\alpha \> 
+ \dots 
\end{eqnarray*}

Checking the coefficients:
\[
 \< \d^i \d^x \>_1 \< \d^j \d^y \d^\mu \> \< \d^\mu \d^z \d^\nu \> \< \d^\nu \d^\alpha \d^\alpha \> : \quad
 \frac{1}{24} c_4 - \frac{1}{24} c_{10} - c_{13} +3 c_{15} +2 c_{18} + c_{20} = 0. 
\]
\[
 \< \d^i \d^\mu \>_1 \< \d^\mu \d^x \d^y \> \< \d^j \d^z \d^\nu \> \< \d^\nu \d^\alpha \d^\alpha \> : \quad
 \frac{1}{24} c_3 - \frac{1}{24} c_{10} + c_{20} = 0. 
\]
\[
 \< \d^\mu \d^\nu \>_1 \< \d^i \d^\mu \d^\nu \> \< \d^j \d^x \d^\alpha \> \< \d^\alpha \d^y \d^z \> : \quad
 \frac{1}{10} c_{10} + c_{16} -2 c_{20} = 0. 
\]
\[
 \< \d^x \d^\mu \>_1 \< \d^\mu \d^i \d^y \> \< \d^j \d^z \d^\nu \> \< \d^\nu \d^\alpha \d^\alpha \> : \quad
 - \frac{1}{12} c_8 +2 c_{13} -6 c_{15} -2 c_{18} = 0. 
\]

\subsection{Checking $\mathfrak{r}_3 (E) =0$}
In the same spirit as the case $l=2$, only the following four
disconnected output graphs will be presented here.

$$ \< \d^x \d^y \d^z \d^i \> \< \d^j \d^\alpha \d^\mu \> \< \d^\mu \d^\alpha \d^\nu \> \< \d^\nu \d^\beta \d^\beta \>, \quad \quad 
\< \d^i \>_1 \< \d^x \d^y \d^\mu \> \< \d^\mu \d^z \d^\nu \> \< \d^\nu \d^i \d^\alpha \> \< \d^\alpha \d^\beta \d^\beta \>, $$
$$ \< \d^i \d^x \d^\mu \d^\mu \> \< \d^y \d^z \d^\nu \> \< \d^\nu \d^j \d^\alpha \> \< \d^\alpha \d^\beta \d^\beta \>, \quad \quad 
\< \d^\mu \>_1 \< \d^\mu \d^i \d^x \> \< \d^y \d^z \d^\nu \> \< \d^\nu \d^j \d^\alpha \> \< \d^\alpha \d^\beta \d^\beta \>. $$
\[
\< \d^x \d^\mu \d^\nu \> \< \d^\mu \d^y \d^z \> \< \d^\nu \>_2 
\mapsto
- \frac{1}{24} \< \d^j \>_1 \< \d^x \d^\mu \d^\nu \> \< \d^\mu \d^y \d^z \> \< \d^i \d^\alpha \d^\nu \> \< \d^\alpha \d^\beta \d^\beta \> 
+ \dots
\]
\begin{eqnarray*} 
& & \< \d^x \d^y \d^z \d^\mu \> \< \d^\mu_1 \>_2 \mapsto  \\ & & \hspace{.5cm}
- \frac{1}{5760}  \< \d^x \d^y \d^z \d^j \> \< \d^i \d^\alpha \d^\mu \> \< \d^\alpha \d^\beta \d^\beta \> \< \d^\mu \d^\nu \d^\nu \> 
+ \frac{1}{960}  \< \d^x \d^y \d^z \d^j \> \< \d^i \d^\alpha \d^\nu \> \< \d^\alpha \d^\beta \d^\nu \> \< \d^\beta \d^\mu \d^\mu \>
+ \dots
\end{eqnarray*}
\begin{eqnarray*} & & 
\< \d^x_1 \d^\mu \>_2 \< \d^\mu \d^y \d^z \> 
\mapsto \\ & & \hspace{.5cm} 
- \frac{1}{24}  \< \d^j \>_1 \< \d^x \d^\mu \d^\beta \> \< \d^\beta \d^i \d^\nu \> \< \d^\nu \d^\alpha \d^\alpha \> \< \d^\mu \d^y \d^z \> 
- \frac{1}{24}  \< \d^j \>_1 \< \d^i \d^\nu \d^\mu \> \< \d^\alpha \d^x \d^\beta \> \< \d^\beta \d^\nu \d^\alpha \> \< \d^\mu \d^y \d^z \> 
\\ & &\hspace{.5cm} 
- \frac{1}{24}  \< \d^x \d^\nu \d^j \> \< \d^\nu \>_1 \< \d^i \d^\alpha \d^\mu \> \< \d^\alpha \d^\beta \d^\beta \> \< \d^\mu \d^y \d^z \> 
- \frac{1}{576}  \< \d^x \d^\nu \d^\nu \d^j \> \< \d^i \d^\alpha \d^\mu \> \< \d^\alpha \d^\beta \d^\beta \> \< \d^\mu \d^y \d^z \> 
+ \dots
\end{eqnarray*}
\begin{eqnarray*} & & 
\< \d^x \d^\mu_1 \>_2 \< \d^\mu \d^y \d^z \> 
\mapsto \\ & & \hspace{.5cm} 
- \frac{1}{24}  \< \d^j \>_1 \< \d^\mu \d^x \d^\beta \> \< \d^\beta \d^i \d^\nu \> \< \d^\nu \d^\alpha \d^\alpha \> \< \d^\mu \d^y \d^z \> 
- \frac{1}{24}  \< \d^j \>_1 \< \d^i \d^\nu \d^x \> \< \d^\alpha \d^\mu \d^\beta \> \< \d^\beta \d^\nu \d^\alpha \> \< \d^\mu \d^y \d^z \> 
\\ & & \hspace{.5cm} 
+ \frac{1}{24}  \< \d^\mu \>_1 \< \d^\mu \d^i \d^x \> \< \d^\nu \d^\beta \d^\alpha \> \< \d^\beta \d^j \d^\alpha \> \< \d^\nu \d^y \d^z \> 
+ \frac{1}{576}  \< \d^i \d^\mu \d^\mu \d^x \> \< \d^\nu \d^\beta \d^\alpha \> \< \d^\beta \d^j \d^\alpha \> \< \d^\nu \d^y \d^z \> 
+ \dots
\end{eqnarray*}
\[
\< \d^x \d^\mu \d^\nu \> \< \d^\mu \>_1 \< \d^\nu \d^y \d^z \>_1 
\mapsto 
\frac{1}{24}  \< \d^\mu \>_1 \< \d^\mu \d^x \d^j \> \< \d^i \d^z \d^\alpha \> \< \d^\alpha \d^\nu \d^y \> \< \d^\nu \d^\beta \d^\beta \> 
+ \dots
\]
\begin{eqnarray*} & & 
\< \d^\mu \>_1 \< \d^\mu \d^x \d^\nu \>_1 \< \d^\nu \d^y \d^z \> 
\mapsto \\ & & \hspace{.5cm}  
\frac{1}{24}  \< \d^j \>_1 \< \d^i \d^\nu \d^\mu \> \< \d^\mu \d^\alpha \d^x \> \< \d^\alpha \d^\beta \d^\beta \> \< \d^\nu \d^y \d^z \> 
- \frac{1}{24}  \< \d^\mu \>_1 \< \d^j \d^x \d^\mu \> \< \d^i \d^\alpha \d^\nu \> \< \d^\alpha \d^\beta \d^\beta \> \< \d^\nu \d^y \d^z \> 
+ \dots
\end{eqnarray*}
\begin{eqnarray*} & & 
\< \d^x \d^y \d^z \d^\mu \d^\mu \d^\nu \> \< \d^\nu \>_1 
\mapsto \\ & & \hspace{.5cm}  
 \< \d^j \>_1 \< \d^i \d^\mu \d^\nu \> \< \d^\nu \d^\alpha \d^\mu \> \< \d^\alpha \d^\beta \d^z \> \< \d^\beta \d^x \d^y \> 
-3  \< \d^\mu \>_1 \< \d^j \d^x \d^\mu \> \< \d^i \d^\alpha \d^\nu \> \< \d^\alpha \d^\beta \d^\nu \> \< \d^\beta \d^y \d^z \> 
+ \dots
\end{eqnarray*}
\[
\< \d^x \d^\mu \d^\mu \d^\nu \> \< \d^\nu \d^y \d^z \>_1 
\mapsto
\frac{1}{24}  \< \d^j \d^x \d^\mu \d^\mu \> \< \d^i \d^z \d^\nu \> \< \d^\nu \d^\alpha \d^y \> \< \d^\alpha \d^\beta \d^\beta \> 
+ \dots
\]
\begin{eqnarray*} & & 
\< \d^x \d^y \d^z \d^\mu \>_1 \< \d^\mu \d^\nu \d^\nu \> 
\mapsto \\ & & \hspace{.5cm} 
- \frac{1}{24}  \< \d^x \d^y \d^z \d^j \> \< \d^i \d^\alpha \d^\mu \> \< \d^\alpha \d^\beta \d^\beta \> \< \d^\mu \d^\nu \d^\nu \> 
-  \< \d^j \>_1 \< \d^i \d^\alpha \d^\mu \> \< \d^\alpha \d^\beta \d^z \> \< \d^\beta \d^x \d^y \> \< \d^\mu \d^\nu \d^\nu \> 
\\ & & 
+3  \< \d^\mu \>_1 \< \d^\mu \d^i \d^x \> \< \d^j \d^\beta \d^\alpha \> \< \d^\beta \d^y \d^z \> \< \d^\alpha \d^\nu \d^\nu \> 
+ \frac{1}{8}  \< \d^i \d^x \d^\mu \d^\mu \> \< \d^j \d^\beta \d^\alpha \> \< \d^\beta \d^y \d^z \> \< \d^\alpha \d^\nu \d^\nu \> 
+ \dots
\end{eqnarray*}
\begin{eqnarray*} & & 
\< \d^x \d^\mu \d^\nu \> \< \d^\mu \d^\nu \d^y \d^z \>_1 
\mapsto \\ & & \hspace{.5cm} 
-  \< \d^j \>_1 \< \d^x \d^\mu \d^\nu \> \< \d^i \d^\alpha \d^z \> \< \d^\alpha \d^\beta \d^y \> \< \d^\beta \d^\mu \d^\nu \> 
+2  \< \d^\mu \>_1 \< \d^\mu \d^i \d^y \> \< \d^x \d^\alpha \d^\nu \> \< \d^j \d^\beta \d^z \> \< \d^\beta \d^\alpha \d^\nu \> 
\\ & & \hspace{.5cm} 
+ \frac{1}{12}  \< \d^i \d^y \d^\mu \d^\mu \> \< \d^x \d^\alpha \d^\nu \> \< \d^j \d^\beta \d^z \> \< \d^\beta \d^\alpha \d^\nu \> 
+ \dots 
\end{eqnarray*}
\[
\< \d^x \d^y \d^z \d^\mu \> \< \d^\mu \d^\nu \d^\nu \>_1 
\mapsto 
\frac{1}{24}  \< \d^x \d^y \d^z \d^j \> \< \d^i \d^\nu \d^\mu \> \< \d^\mu \d^\alpha \d^\nu \> \< \d^\alpha \d^\beta \d^\beta \> 
+ \dots
\]
\begin{eqnarray*} & & 
\< \d^x \d^y \d^\mu \> \< \d^\mu \d^\nu \d^\nu \d^z \>_1 
\mapsto \\ & & \hspace{.5cm} 
- \frac{1}{24}  \< \d^x \d^y \d^\mu \> \< \d^i \d^\alpha \d^\mu \> \< \d^\alpha \d^\beta \d^\beta \> \< \d^j \d^\nu \d^\nu \d^z \> 
-  \< \d^x \d^y \d^\mu \> \< \d^j \>_1 \< \d^i \d^\alpha \d^z \> \< \d^\alpha \d^\beta \d^\nu \> \< \d^\beta \d^\mu \d^\nu \> 
\\ & & \hspace{.5cm} 
+  \< \d^\mu \>_1 \< \d^\mu \d^i \d^z \> \< \d^j \d^\beta \d^\nu \> \< \d^\beta \d^\nu \d^\alpha \> \< \d^x \d^y \d^\alpha \> 
+ \frac{1}{24}  \< \d^i \d^z \d^\mu \d^\mu \> \< \d^j \d^\beta \d^\nu \> \< \d^\beta \d^\nu \d^\alpha \> \< \d^x \d^y \d^\alpha \> 
+ \dots
\end{eqnarray*}

Checking the coefficients:
\[
 \< \d^x \d^y \d^z \d^i \> \< \d^j \d^\alpha \d^\mu \> \< \d^\mu \d^\alpha \d^\nu \> \< \d^\nu \d^\beta \d^\beta \> : \quad
 \frac{1}{1152} c_2 - \frac{1}{24} c_{15} + \frac{1}{24} c_{19} = 0.
\]
\[
 \< \d^i \>_1 \< \d^x \d^y \d^\mu \> \< \d^\mu \d^z \d^\nu \> \< \d^\nu \d^i \d^\alpha \> \< \d^\alpha \d^\beta \d^\beta \> : \quad
 - \frac{1}{24} c_1 - \frac{1}{12} c_3 - \frac{1}{12} c_4 + \frac{1}{24} c_{11} + c_{12} - c_{15} - c_{18} - c_{20} = 0.
\]
\[
 \< \d^i \d^x \d^\mu \d^\mu \> \< \d^y \d^z \d^\nu \> \< \d^\nu \d^j \d^\alpha \> \< \d^\alpha \d^\beta \d^\beta \> : \quad
 - \frac{1}{576} c_3 + \frac{1}{576} c_4 + \frac{1}{24} c_{14} + \frac{1}{8} c_{15} + \frac{1}{12} c_{18} = 0. 
\]
\[
 \< \d^\mu \>_1 \< \d^\mu \d^i \d^x \> \< \d^y \d^z \d^\nu \> \< \d^\nu \d^j \d^\alpha \> \< \d^\alpha \d^\beta \d^\beta \> : \quad
 - \frac{1}{24} c_3 + \frac{1}{24} c_4 + \frac{1}{24} c_7 - \frac{1}{24} c_{11} -3 c_{12} +3 c_{15} +2 c_{18} + c_{20} = 0. 
\]

\subsection{Checking $\mathfrak{r}_4(E)=0$}
Since $l=4$ case is new, the calculation is presented.

\begin{eqnarray*} & & 
\< \d^x \d^\mu \d^\nu \> \< \d^\mu \d^y \d^z \> \< \d^\nu \>_2 
\mapsto \\ & & \hspace{.5cm} 
- \frac{1}{5760}  \< \d^x \d^\mu \d^j \> \< \d^\mu \d^y \d^z \> \< \d^i \d^\beta \d^\nu \> \< \d^\beta \d^\gamma \d^\gamma \> \< \d^\nu \d^\alpha \d^\alpha \> 
\\ & & \hspace{.5cm} 
+ \frac{1}{960}  \< \d^x \d^\mu \d^j \> \< \d^\mu \d^y \d^z \> \< \d^i \d^\beta \d^\alpha \> \< \d^\beta \d^\gamma \d^\alpha \> \< \d^\gamma \d^\nu \d^\nu \> 
\\ & & \hspace{.5cm} 
+ \frac{1}{576}  \< \d^x \d^\mu \d^\nu \> \< \d^\mu \d^y \d^z \> \< \d^i \d^\beta \d^\nu \> \< \d^\beta \d^\gamma \d^\gamma \> \< \d^j \d^\alpha \d^\alpha \> 
\end{eqnarray*}
\begin{eqnarray*} & & 
\< \d^x_1 \d^\mu \>_2 \< \d^\mu \d^y \d^z \> 
\mapsto \\ & & \hspace{.5cm} 
- \frac{1}{1920}  \< \d^i \d^\nu \d^\alpha \> \< \d^\alpha \d^\beta \d^x \> \< \d^\beta \d^\gamma \d^\gamma \> \< \d^\mu \d^\mu \d^\nu \> \< \d^j \d^y \d^z \> 
\\ & & \hspace{.5cm} 
+ \frac{1}{720}  \< \d^x \d^\mu \d^\nu \> \< \d^i \d^\nu \d^\alpha \> \< \d^\alpha \d^\beta \d^\mu \> \< \d^\beta \d^\gamma \d^\gamma \> \< \d^j \d^y \d^z \> 
\\ & & \hspace{.5cm} 
+ \frac{1}{576}  \< \d^i \d^\nu \d^\alpha \> \< \d^\alpha \d^\beta \d^\nu \> \< \d^\beta \d^\gamma \d^\mu \> \< \d^\gamma \d^x \d^\mu \> \< \d^j \d^y \d^z \> 
\\ & & \hspace{.5cm} 
+ \frac{1}{576}  \< \d^x \d^\mu \d^\alpha \> \< \d^\alpha \d^i \d^\beta \> \< \d^\beta \d^\gamma \d^\gamma \> \< \d^j \d^\nu \d^\nu \> \< \d^\mu \d^y \d^z \> 
\\ & & \hspace{.5cm} 
+ \frac{1}{576}  \< \d^i \d^\beta \d^\mu \> \< \d^\gamma \d^x \d^\alpha \> \< \d^\alpha \d^\beta \d^\gamma \> \< \d^j \d^\nu \d^\nu \> \< \d^\mu \d^y \d^z \> 
\\ & & \hspace{.5cm} 
+ \frac{1}{576}  \< \d^i \d^\beta \d^\mu \> \< \d^\beta \d^\gamma \d^\gamma \> \< \d^x \d^\alpha \d^\nu \> \< \d^\alpha \d^j \d^\nu \> \< \d^\mu \d^y \d^z \> 
\end{eqnarray*}
\begin{eqnarray*} & & 
\< \d^x \d^\mu_1 \>_2 \< \d^\mu \d^y \d^z \> 
\mapsto \\ & & \hspace{.5cm} 
- \frac{1}{5760}  \< \d^i \d^x \d^\alpha \> \< \d^\alpha \d^\beta \d^\mu \> \< \d^\beta \d^\gamma \d^\gamma \> \< \d^\mu \d^\nu \d^\nu \> \< \d^j \d^y \d^z \> 
\\ & & \hspace{.5cm} 
+ \frac{1}{960}  \< \d^i \d^x \d^\alpha \> \< \d^\alpha \d^\beta \d^\nu \> \< \d^\beta \d^\gamma \d^\nu \> \< \d^\gamma \d^\mu \d^\mu \> \< \d^j \d^y \d^z \> 
\\ & & \hspace{.5cm} 
+ \frac{1}{576}  \< \d^\mu \d^x \d^\alpha \> \< \d^\alpha \d^i \d^\beta \> \< \d^\beta \d^\gamma \d^\gamma \> \< \d^j \d^\nu \d^\nu \> \< \d^\mu \d^y \d^z \> 
\\ & & \hspace{.5cm} 
+ \frac{1}{576}  \< \d^i \d^\beta \d^x \> \< \d^\gamma \d^\mu \d^\alpha \> \< \d^\alpha \d^\beta \d^\gamma \> \< \d^j \d^\nu \d^\nu \> \< \d^\mu \d^y \d^z \> 
\\ & & \hspace{.5cm} 
+ \frac{1}{576}  \< \d^i \d^\beta \d^x \> \< \d^\beta \d^\gamma \d^\gamma \> \< \d^\mu \d^\alpha \d^\nu \> \< \d^\alpha \d^j \d^\nu \> \< \d^\mu \d^y \d^z \> 
\end{eqnarray*}
\begin{eqnarray*} & & 
\< \d^x \d^y \d^z \d^\mu \>_1 \< \d^\mu \d^\nu \d^\nu \> 
\mapsto \\ & & \hspace{.5cm} 
+ \frac{1}{24}  \< \d^i \d^\alpha \d^\nu \> \< \d^\nu \d^\beta \d^\alpha \> \< \d^\beta \d^\gamma \d^z \> \< \d^\gamma \d^x \d^y \> \< \d^j \d^\mu \d^\mu \> 
\\ & & \hspace{.5cm} 
- \frac{1}{8}  \< \d^i \d^y \d^\alpha \> \< \d^\alpha \d^\beta \d^x \> \< \d^\beta \d^\gamma \d^\gamma \> \< \d^j \d^z \d^\mu \> \< \d^\mu \d^\nu \d^\nu \> 
\\ & & \hspace{.5cm} 
- \frac{1}{8}  \< \d^i \d^\mu \d^\alpha \> \< \d^\alpha \d^\beta \d^x \> \< \d^\beta \d^\gamma \d^\gamma \> \< \d^j \d^y \d^z \> \< \d^\mu \d^\nu \d^\nu \> 
\\ & & \hspace{.5cm} 
+ \frac{1}{8}  \< \d^i \d^\beta \d^x \> \< \d^\beta \d^\gamma \d^\gamma \> \< \d^j \d^\alpha \d^\mu \> \< \d^\alpha \d^y \d^z \> \< \d^\mu \d^\nu \d^\nu \> 
\\ & & \hspace{.5cm} 
+ \frac{1}{24}  \< \d^i \d^\beta \d^\mu \> \< \d^\beta \d^\gamma \d^\gamma \> \< \d^j \d^\alpha \d^z \> \< \d^\alpha \d^x \d^y \> \< \d^\mu \d^\nu \d^\nu \> 
\\ & & \hspace{.5cm} 
+ \frac{1}{24}  \< \d^i \d^\beta \d^\mu \> \< \d^\beta \d^\gamma \d^z \> \< \d^\gamma \d^x \d^y \> \< \d^j \d^\alpha \d^\alpha \> \< \d^\mu \d^\nu \d^\nu \> 
\end{eqnarray*}
\begin{eqnarray*} & & 
\< \d^x \d^\mu \d^\nu \> \< \d^\mu \d^\nu \d^y \d^z \>_1 
\mapsto \\ & & \hspace{.5cm} 
- \frac{1}{24}  \< \d^x \d^\mu \d^\nu \> \< \d^i \d^\nu \d^\alpha \> \< \d^\alpha \d^\beta \d^\mu \> \< \d^\beta \d^\gamma \d^\gamma \> \< \d^j \d^y \d^z \> 
\\ & & \hspace{.5cm} 
- \frac{1}{24}  \< \d^x \d^\mu \d^\nu \> \< \d^i \d^x \d^\alpha \> \< \d^\alpha \d^\beta \d^y \> \< \d^\beta \d^\gamma \d^\gamma \> \< \d^j \d^\mu \d^\nu \> 
\\ & & \hspace{.5cm} 
+ \frac{1}{12}  \< \d^x \d^\mu \d^\nu \> \< \d^i \d^\beta \d^y \> \< \d^\beta \d^\gamma \d^\gamma \>\< \d^j \d^\alpha \d^z \> \< \d^\alpha \d^\mu \d^\nu \>  
\\ & & \hspace{.5cm} 
+ \frac{1}{24}  \< \d^x \d^\mu \d^\nu \> \< \d^i \d^\beta \d^z \> \< \d^\beta \d^\gamma \d^y \> \< \d^\gamma \d^\mu \d^\nu \> \< \d^j \d^\alpha \d^\alpha \> 
\end{eqnarray*}
\begin{eqnarray*} & & 
\< \d^x \d^y \d^\mu \> \< \d^\mu \d^\nu \d^\nu \d^z \>_1 
\mapsto \\ & & \hspace{.5cm} 
+ \frac{1}{24}  \< \d^x \d^y \d^j \> \< \d^i \d^\alpha \d^\nu \> \< \d^\nu \d^\beta \d^\alpha \> \< \d^\beta \d^\gamma \d^z \> \< \d^\gamma \d^\mu \d^\mu \> 
\\ & & \hspace{.5cm} 
- \frac{1}{24}  \< \d^x \d^y \d^\mu \> \< \d^i \d^z \d^\alpha \> \< \d^\alpha \d^\beta \d^\mu \> \< \d^\beta \d^\gamma \d^\gamma \> \< \d^j \d^\nu \d^\nu \> 
\\ & & \hspace{.5cm} 
- \frac{1}{24}  \< \d^x \d^y \d^\mu \> \< \d^i \d^\nu \d^\alpha \> \< \d^\alpha \d^\beta \d^\nu \> \< \d^\beta \d^\gamma \d^\gamma \> \< \d^j \d^\mu \d^z \> 
\\ & & \hspace{.5cm} 
+ \frac{1}{24}  \< \d^x \d^y \d^\mu \> \< \d^i \d^\beta \d^\mu \> \< \d^\beta \d^\gamma \d^\gamma \> \< \d^j \d^\alpha \d^z \> \< \d^\alpha \d^\nu \d^\nu \> 
\\ & & \hspace{.5cm} 
+ \frac{1}{24}  \< \d^x \d^y \d^\mu \> \< \d^i \d^\beta \d^z \> \< \d^\beta \d^\gamma \d^\gamma \> \< \d^j \d^\alpha \d^\nu \> \< \d^\alpha \d^\mu \d^\nu \> 
\\ & & \hspace{.5cm} 
+ \frac{1}{24}  \< \d^x \d^y \d^\mu \> \< \d^i \d^\beta \d^z \> \< \d^\beta \d^\gamma \d^\nu \> \< \d^\gamma \d^\mu \d^\nu \> \< \d^j \d^\alpha \d^\alpha \>
\end{eqnarray*}
The other graphs all have $\mathfrak{r}_4(\Gamma)=0$.

Now let's check all the coefficients:
\[
 \< \d^x \d^\mu \d^\nu \> \< \d^\mu \d^y \d^z \> \< \d^i \d^\beta \d^\nu \> \< \d^\beta \d^\gamma \d^\gamma \> \< \d^j \d^\alpha \d^\alpha \>: \quad
 \frac{1}{576} c_1 + \frac{1}{288} c_3 + \frac{1}{288} c_4 + \frac{1}{12} c_{15} + \frac{1}{24} c_{18} = 0.
\]
\[
 \< \d^i \d^\nu \d^\alpha \> \< \d^\alpha \d^\beta \d^x \> \< \d^\beta \d^\gamma \d^\gamma \> \< \d^\mu \d^\mu \d^\nu \> \< \d^j \d^y \d^z \> : \quad
 \frac{1}{384} c_3 + \frac{1}{1152} c_4 - \frac{1}{8} c_{15} - \frac{1}{24} c_{18} + \frac{1}{24} c_{20} = 0.
\]
\[
 \< \d^x \d^\mu \d^j \> \< \d^\mu \d^y \d^z \> \< \d^i \d^\beta \d^\nu \> \< \d^\beta \d^\gamma \d^\gamma \> \< \d^\nu \d^\alpha \d^\alpha \> : \quad
 \frac{1}{1152} c_1 + \frac{1}{24} c_{15} - \frac{1}{24} c_{20} = 0. 
\]
\[
 \< \d^j \d^x \d^\mu \> \< \d^\mu \d^\nu \d^\nu \> \< \d^i \d^y \d^\alpha \> \< \d^\alpha \d^\beta \d^z \> \< \d^\beta \d^\gamma \d^\gamma \> : \quad
 \frac{1}{576} c_3 - \frac{1}{576} c_4 - \frac{1}{4} c_{15} - \frac{1}{8} c_{18} = 0. 
\]

\subsection{Conclusion}
By Lemma~1 in \cite{ypL2}, $\mathfrak{r}_l(E)=0$ for $l \ge 3$ 
(respectively $4,5$) for $(g,n,k) = (2,1,2)$, 
(respectively $(2,2,2), (2,3,2)$).
Therefore, Conjectures~1 and 2 are proved. 
By a Betti number calculation of E.~Getzler \cite{eG2}, 
they are the only tautological equations for the corresponding $(g,n,k)$. 
Therefore, Conjecture~3 also holds.


\begin{thebibliography}{Co2}

\bibitem{AL}
D.~Arcara, Y.-P.~Lee, \emph{Tautological equation in $\Mbar_{3,1}$ via
invariance conjecture}, math.AG/0503184.

\bibitem{BP}
P.~Belorousski, R.~Pandharipande, \emph{A descendent relation in genus 2},
Ann. Scuola Norm. Sup. Pisa Cl. Sci. (4) 29 (2000), no. 1, 171--191.

\bibitem{eG1}
E.~Getzler, \emph{Intersection theory on $\overline{\mathcal{M}}_{1,4}$ and
elliptic Gromov-Witten invariants}, J. Amer. Math. Soc. 10 (1997), no. 4,
973--998.

\bibitem{eG2}
E.~Getzler, \emph{Topological recursion relations in genus $2$},
Integrable systems and algebraic geometry (Kobe/Kyoto, 1997),  73--106,
World Sci. Publishing, River Edge, NJ, 1998.

\bibitem{GL}
A.~Givental, Y.-P.~Lee, preliminary draft. 

\bibitem{ypL1}
Y.-P.~Lee, \emph{Witten's conjecture and Virasoro conjecture up to genus two},
math.AG/0310442. To appear in the proceedings of the conference 
"Gromov-Witten Theory of Spin Curves and Orbifolds", Contemp. Math., AMS.

\bibitem{ypL2}
Y.-P.~Lee, \emph{Invariance of tautological equations I: conjectures and 
applications},math.AG/0604318. 

\bibitem{ypL3}
Y.-P.~Lee, \emph{Invariance of tautological equations II: Gromov--Witten 
theory}, preprint available at http://www.math.utah.edu/\~{}yplee/research/.

\bibitem{ypL4}
Y.-P.~Lee, \emph{Witten's conjecture, Virasoro conjecture, and invariance 
of tautological relations}, math.AG/0311100.


\end{thebibliography}
\end{document}